\documentclass{gtart}
\def\ifplaintex{\expandafter\ifx\csname documentclass\endcsname\relax}

\def\gtp{{\mathsurround=0pt\it $\cal G\mskip-2mu$eometry \&\ 
$\cal T\!\!$opology $\cal P\!$ublications}}  % GT publications

\def\recd{{\small Received:\qua\receiveddate\ifx\reviseddate\relax
\else\qquad Revised:\qua\reviseddate\fi\par}} 

\def\lognumber#1{\def\thelognumber{#1}}
\def\volumenumber#1{\def\thevolumenumber{#1}}
\def\volumeyear#1{\def\thevolumeyear{#1}}
\def\papernumber#1{\def\thepapernumber{#1}}
\def\pagenumbers#1#2{\def\startpage{#1}\def\finishpage{#2}}
\def\published#1{\def\publishdate{#1}}

\def\received#1{\def\receiveddate{#1}}
\def\revised#1{\def\reviseddate{#1}}
\def\accepted#1{\def\accepteddate{#1}}
\def\asciititle#1{\def\theasciititle{#1}}
\def\covertitle#1{\def\thecovertitle{#1}}

\def\asciiaddress#1{\def\theasciiaddress{#1}}

\long\def\asciiabstract#1{\long\def\theasciiabstract{#1}}

\let\\\par\let\thelognumber\relax\let\thevolumenumber\relax
\let\thepapernumber\relax\let\thevolumeyear\relax\let\startpage\relax
\let\finishpage\relax\let\publishdate\relax\let\receiveddate\relax
\let\reviseddate\relax\let\accepteddate\relax\let\theasciititle\relax
\let\thecovertitle\relax\let\theasciiauthors\relax\let\theasciiaddress\relax
\let\theasciiabstract\relax

\let\theasciiemail\relax

\ifplaintex
\font\logobig=cmssbx10 scaled 3836
\font\logomed=cmssbx10 scaled 2557
\else
\font\logobig=cmssbx10 scaled 4200
\font\logomed=cmssbx10 scaled 2800
\fi

\long\def\makeagttitle{   %%% start of definition of \makeagttitle
\count0=\startpage
\agt\hfill      %   Journal title (top left) 
\hbox to 45truept{\vbox to 0pt{\vglue -13truept{\logomed A\kern -.37em{\logobig 
T}\kern -.38em G}\vss}\hss}
\break
{\small Volume \thevolumenumber\ (\thevolumeyear)
\startpage--\finishpage\nl
Published: \publishdate}

\vglue .25truein

{\parskip=0pt\leftskip 0pt plus
1fil\def\\{\par\smallskip}{\Large\bf\thetitle}\par\medskip} \vglue
0.05truein

{\parskip=0pt\leftskip 0pt plus 1fil\def\\{\par}{\sc\theauthors}
\par\medskip}%
 
\vglue 0.03truein 

{\small\leftskip 25truept\rightskip 25truept{\bf Abstract}\stdspace\theabstract

{\bf AMS Classification}\stdspace\theprimaryclass
\ifx\thesecondaryclass\relax\else; \thesecondaryclass\fi\par
{\bf Keywords}\stdspace \thekeywords\par}\vglue 7truept

}   %%%% end of definition of \makeagttitle

\ifplaintex
\font\phead=cmsl9 scaled 950
\font\pnum=cmbx10 scaled 913
\font\pfoot=cmsl9 scaled 950
\headline{\vbox to 0pt{\vskip -4.5mm\line{\small\phead\ifnum
\count0=\startpage ISSN 1472-2739 (on-line) 1472-2747 (printed)
\hfill {\pnum\folio}\else\ifodd\count0\def\\{ }% 
\ifx\theshorttitle\relax\thetitle\else\theshorttitle\fi\hfill{\pnum\folio}
\else\def\\{ and }{\pnum\folio}\hfill\ifx\theshortauthors\relax\theauthors
\else\theshortauthors\fi\fi\fi}\vss}}
\footline{\vbox to 0pt{\vglue 0mm\line{\small\pfoot\ifnum\count0=\startpage
\copyright\ \gtp\hfill\else
\agt, Volume \thevolumenumber\ (\thevolumeyear)\hfill\fi}\vss}}
\else
\font=cmsl9 scaled 1050
\font\lnum=cmbx10 
\font=cmsl9 scaled 1050
\makeatletter
\def\@oddhead{{\small\lhead\ifnum\count0=\startpage ISSN 1472-2739 
(on-line) 1472-2747 (printed)\hfill {\lnum\number\count0}\else\ifodd\count0
\def\\{ }\ifx\theshorttitle\relax \thetitle \else\theshorttitle\fi\hfill
{\lnum\number\count0}\else\def\\{ and }{\lnum\number\count0}
\hfill\ifx\theshortauthors\relax 
\theauthors\else\theshortauthors\fi\fi\fi}}\def\@evenhead{\@oddhead}
\def\@oddfoot{\small\lfoot\ifnum\count0=\startpage\copyright\ \gtp\hfill\else
\agt, Volume \thevolumenumber\ (\thevolumeyear)\hfill\fi}
\def\@evenfoot{\@oddfoot}
\makeatother
\fi
\let\maketitlepage\makeagttitle
\let\makeshorttitle\maketitlepage
\let\maketitle\maketitlepage

\newwrite\gtoutfile
\long\gdef\makeheadfile{  %%% start of definition of \makeheadfile
{\def\\{, }\def\s{ }
\immediate\openout\gtoutfile head.xxx
\immediate\write\gtoutfile{To: math@arxiv.org}
\immediate\write\gtoutfile{Subject: put OR rep NNNNN:ppppp}
\immediate\write\gtoutfile{--text follows this line--}
\immediate\write\gtoutfile{Proxy-for: \ifx\theasciiauthors\relax
\theauthors\else\theasciiauthors\fi\s<\ifx\theasciiemail\relax\theemail\else\theasciiemail\fi>}
\immediate\write\gtoutfile{\noexpand\\}
\immediate\write\gtoutfile{Authors: \ifx\theasciiauthors\relax
\theauthors\else\theasciiauthors\fi}
{\def\\{ }\immediate\write\gtoutfile{Title: \ifx\theasciititle\relax
\thetitle\else\theasciititle\fi}}
\immediate\write\gtoutfile{Subj-class: GT or SG, GR etc}
\immediate\write\gtoutfile{MSC-class: \theprimaryclass\ifx\thesecondaryclass\relax\else, \thesecondaryclass\fi}
\immediate\write\gtoutfile{Journal-ref: Algebr. Geom. Topol. \thevolumenumber\s
(\thevolumeyear) \startpage-\finishpage}
\immediate\write\gtoutfile{Comments: Published by Algebraic and
Geometric Topology at}
\immediate\write\gtoutfile{\s\s\s  http://www.maths.warwick.ac.uk/agt/AGTVol\thevolumenumber/agt-\thevolumenumber-\thepapernumber.abs.html}
\immediate\write\gtoutfile{\noexpand\\}
\immediate\write\gtoutfile{}
\ifx\theasciiabstract\relax
\immediate\write\gtoutfile{\theabstract}\else
\immediate\write\gtoutfile{\theasciiabstract}\fi
\immediate\write\gtoutfile{}
\immediate\write\gtoutfile{\noexpand\\}
\immediate\write\gtoutfile{}
\immediate\closeout\gtoutfile}}  %%% end of definition of \makeheadfile

\def\maketitlepage{\makeagttitle\makeheadfile}
\let\makeshorttitle\maketitlepage
\let\maketitle\maketitlepage

\def\ifplaintex{\expandafter\ifx\csname documentclass\endcsname\relax}

\def\gtp{{\mathsurround=0pt\it $\cal G\mskip-2mu$eometry \&\ 
$\cal T\!\!$opology $\cal P\!$ublications}}  % GT publications

\def\recd{{\small Received:\qua\receiveddate\ifx\reviseddate\relax
\else\qquad Revised:\qua\reviseddate\fi\par}} 

\def\lognumber#1{\def\thelognumber{#1}}
\def\volumenumber#1{\def\thevolumenumber{#1}}
\def\volumeyear#1{\def\thevolumeyear{#1}}
\def\papernumber#1{\def\thepapernumber{#1}}
\def\pagenumbers#1#2{\def\startpage{#1}\def\finishpage{#2}}
\def\published#1{\def\publishdate{#1}}

\def\received#1{\def\receiveddate{#1}}
\def\revised#1{\def\reviseddate{#1}}
\def\accepted#1{\def\accepteddate{#1}}
\def\asciititle#1{\def\theasciititle{#1}}
\def\covertitle#1{\def\thecovertitle{#1}}

\def\asciiaddress#1{\def\theasciiaddress{#1}}

\long\def\asciiabstract#1{\long\def\theasciiabstract{#1}}

\let\\\par\let\thelognumber\relax\let\thevolumenumber\relax
\let\thepapernumber\relax\let\thevolumeyear\relax\let\startpage\relax
\let\finishpage\relax\let\publishdate\relax\let\receiveddate\relax
\let\reviseddate\relax\let\accepteddate\relax\let\theasciititle\relax
\let\thecovertitle\relax\let\theasciiauthors\relax\let\theasciiaddress\relax
\let\theasciiabstract\relax

\let\theasciiemail\relax

\ifplaintex
\font\logobig=cmssbx10 scaled 3836
\font\logomed=cmssbx10 scaled 2557
\else
\font\logobig=cmssbx10 scaled 4200
\font\logomed=cmssbx10 scaled 2800
\fi

\long\def\makeagttitle{   %%% start of definition of \makeagttitle
\count0=\startpage
\agt\hfill      %   Journal title (top left) 
\hbox to 45truept{\vbox to 0pt{\vglue -13truept{\logomed A\kern -.37em{\logobig 
T}\kern -.38em G}\vss}\hss}
\break
{\small Volume \thevolumenumber\ (\thevolumeyear)
\startpage--\finishpage\nl
Published: \publishdate}

\vglue .25truein

{\parskip=0pt\leftskip 0pt plus
1fil\def\\{\par\smallskip}{\Large\bf\thetitle}\par\medskip} \vglue
0.05truein

{\parskip=0pt\leftskip 0pt plus 1fil\def\\{\par}{\sc\theauthors}
\par\medskip}%
 
\vglue 0.03truein 

{\small\leftskip 25truept\rightskip 25truept{\bf Abstract}\stdspace\theabstract

{\bf AMS Classification}\stdspace\theprimaryclass
\ifx\thesecondaryclass\relax\else; \thesecondaryclass\fi\par
{\bf Keywords}\stdspace \thekeywords\par}\vglue 7truept

}   %%%% end of definition of \makeagttitle

\ifplaintex
\font\phead=cmsl9 scaled 950
\font\pnum=cmbx10 scaled 913
\font\pfoot=cmsl9 scaled 950
\headline{\vbox to 0pt{\vskip -4.5mm\line{\small\phead\ifnum
\count0=\startpage ISSN 1472-2739 (on-line) 1472-2747 (printed)
\hfill {\pnum\folio}\else\ifodd\count0\def\\{ }% 
\ifx\theshorttitle\relax\thetitle\else\theshorttitle\fi\hfill{\pnum\folio}
\else\def\\{ and }{\pnum\folio}\hfill\ifx\theshortauthors\relax\theauthors
\else\theshortauthors\fi\fi\fi}\vss}}
\footline{\vbox to 0pt{\vglue 0mm\line{\small\pfoot\ifnum\count0=\startpage
\copyright\ \gtp\hfill\else
\agt, Volume \thevolumenumber\ (\thevolumeyear)\hfill\fi}\vss}}
\else
\font=cmsl9 scaled 1050
\font\lnum=cmbx10 
\font=cmsl9 scaled 1050
\makeatletter
\def\@oddhead{{\small\lhead\ifnum\count0=\startpage ISSN 1472-2739 
(on-line) 1472-2747 (printed)\hfill {\lnum\number\count0}\else\ifodd\count0
\def\\{ }\ifx\theshorttitle\relax \thetitle \else\theshorttitle\fi\hfill
{\lnum\number\count0}\else\def\\{ and }{\lnum\number\count0}
\hfill\ifx\theshortauthors\relax 
\theauthors\else\theshortauthors\fi\fi\fi}}\def\@evenhead{\@oddhead}
\def\@oddfoot{\small\lfoot\ifnum\count0=\startpage\copyright\ \gtp\hfill\else
\agt, Volume \thevolumenumber\ (\thevolumeyear)\hfill\fi}
\def\@evenfoot{\@oddfoot}
\makeatother
\fi
\let\maketitlepage\makeagttitle
\let\makeshorttitle\maketitlepage
\let\maketitle\maketitlepage

\newwrite\gtoutfile
\long\gdef\makeheadfile{  %%% start of definition of \makeheadfile
{\def\\{, }\def\s{ }
\immediate\openout\gtoutfile head.xxx
\immediate\write\gtoutfile{To: math@arxiv.org}
\immediate\write\gtoutfile{Subject: put OR rep NNNNN:ppppp}
\immediate\write\gtoutfile{--text follows this line--}
\immediate\write\gtoutfile{Proxy-for: \ifx\theasciiauthors\relax
\theauthors\else\theasciiauthors\fi\s<\ifx\theasciiemail\relax\theemail\else\theasciiemail\fi>}
\immediate\write\gtoutfile{\noexpand\\}
\immediate\write\gtoutfile{Authors: \ifx\theasciiauthors\relax
\theauthors\else\theasciiauthors\fi}
{\def\\{ }\immediate\write\gtoutfile{Title: \ifx\theasciititle\relax
\thetitle\else\theasciititle\fi}}
\immediate\write\gtoutfile{Subj-class: GT or SG, GR etc}
\immediate\write\gtoutfile{MSC-class: \theprimaryclass\ifx\thesecondaryclass\relax\else, \thesecondaryclass\fi}
\immediate\write\gtoutfile{Journal-ref: Algebr. Geom. Topol. \thevolumenumber\s
(\thevolumeyear) \startpage-\finishpage}
\immediate\write\gtoutfile{Comments: Published by Algebraic and
Geometric Topology at}
\immediate\write\gtoutfile{\s\s\s  http://www.maths.warwick.ac.uk/agt/AGTVol\thevolumenumber/agt-\thevolumenumber-\thepapernumber.abs.html}
\immediate\write\gtoutfile{\noexpand\\}
\immediate\write\gtoutfile{}
\ifx\theasciiabstract\relax
\immediate\write\gtoutfile{\theabstract}\else
\immediate\write\gtoutfile{\theasciiabstract}\fi
\immediate\write\gtoutfile{}
\immediate\write\gtoutfile{\noexpand\\}
\immediate\write\gtoutfile{}
\immediate\closeout\gtoutfile}}  %%% end of definition of \makeheadfile

\def\maketitlepage{\makeagttitle\makeheadfile}
\let\makeshorttitle\maketitlepage
\let\maketitle\maketitlepage

\lognumber{27}
\volumenumber{1}
\volumeyear{2001}
\papernumber{27}
\pagenumbers{519}{548}
\received{9 October 2000}
\revised{4 July 2001}
\accepted{30 September 2001}
\published{5 Octoberber 2001}

\begin{document}
\def\Nil{{\cal N}il}
\def\ra{\rightarrow}
\def\Sq{{\rm Sq}}
\newtheorem{thm}{Th\'eor\`eme}[section]
\newtheorem{conj}{Conjecture}[section]
\newtheorem{cor}[thm]{Corollaire}
\newtheorem{lem}[thm]{Lemme}
\newtheorem{prop}[thm]{Proposition}
\theoremstyle{remark}
\newtheorem{defn}[thm]{D\'efinition}
\newtheorem{rem}[thm]{Remarque}
\def\F{{\bf F}_2}
\def\vfl{\downarrow}
\def\diagram#1{\def\normalbaselines{\baselineskip=0pt
\lineskip=10pt\lineskiplimit=1pt} \matrix {#1}}
\def\un{${\cal U}/{\cal N}il$}
\def\U{{\cal U}}
\def\lam{\lambda}

\title{La filtration de Krull de la cat\' egorie ${\cal U}$\\et la
cohomologie des espaces}
\covertitle{La filtration de Krull de la cat\noexpand\'egorie 
${\noexpand\cal U}$\\et 
la cohomologie des espaces}
\asciititle{La filtration de Krull de la categorie  U et la
cohomologie des espaces}
  \author{Lionel Schwartz}
  \address {Universit\'e Paris-Nord, Institut Galil\'ee, LAGA,
UMR 7539 du CNRS\\Av. J.-B. Cl\'ement, 93430, Villetaneuse,
France}
  \asciiaddress {Universite Paris-Nord, Institut Galilee, LAGA,
UMR 7539 du CNRS\\Av. J.-B. Clement, 93430, Villetaneuse,
France}
  \email {schwartz@math.univ-paris13.fr}

  \begin{abstract}
  Cet article d\'emontre une variante d'une conjecture due
\`a N. Kuhn. Cette conjecture s'exprime \`a l'aide de la
filtration de Krull de la cat\' egorie ${\cal U}$ des modules
instables. Notons ${\cal U}_n$, $n \geq 0$, le $n$-i\`eme terme de
cette filtration. La cat\' egorie ${\cal U}$ est la plus petite
sous cat\'egorie \'epaisse contenant les cat\'egories ${\cal U}_n$
et stable par colimite [7]. La cat\'egorie ${\cal U}_0$ est celle
des modules localement finis, c'est-\`a-dire limite directe de
modules finis.

On entend par  sous cat\'egorie \'epaisse une sous-cat\'egorie stable par
sous-objet et quotient, et telle que pour toute suite exacte
courte,
si le premier et le troisi\`eme terme sont dans la sous-cat\'egorie, alors
le terme central l'est aussi.

La conjecture s'\'enonce comme suit, soit $X$ un espace, alors~:

$\bullet$\qua soit $H^*X \in {\cal U}_0$,

$\bullet$\qua soit $H^*X \not \in {\cal U}_n$, pour tout $n$.

Par exemple la cohomologie d'un espace de dimension finie, ou celle de son espace
des lacets sont toujours dans ${\cal U}_0$. Alors
que la cohomologie du classifiant d'un groupe fini, d'ordre
divisible par $2$ n'est, elle, dans aucune des sous-cat\'egories ${\cal U}_n$.

On d\'emontre cette conjecture, modulo l'hypoth\`ese
suppl\'ementaire que tous les quotients
de la filtration nilpotente ont un nombre fini
de g\'en\'erateurs. Cette condition implique en particulier que la cohomologie  est
de dimension finie en chaque degr\'e. Mais elle est plus forte, et
 assure les conditions
d'application du th\'eor\`eme de Lannes sur la cohomologie des
espaces fonctionnels. Ce th\'eor\`eme est n\'ecessaire pour
appliquer la r\'eduction de Kuhn [3].

Par commodit\'e on ne consid\`erera dans cet article que le cas $p=2$.
\end{abstract}

\asciiabstract{This paper proves a particular case of a conjecture of N. Kuhn. This
conjecture is as follows. Consider the Gabriel-Krull filtration of the
category U of unstable modules.

Let U_n, n>=0, be the n-th step of
this filtration. The category U is the smallest thick sub-category
that contains all sub-categories U_n and is stable under colimit
[L. Schwartz, Unstable modules over the Steenrod algebra and
Sullivan's fixed point set conjecture, Chicago Lectures in Mathematics
Series (1994)]. The category U_0 is the one of locally finite modules,
i.e. the modules that are direct limit of finite modules.  The
conjecture is as follows, let X be a space then : 

* either H^*X is locally finite,
* or H^*X does not belong to U_n, for all n.

As an example the cohomology of a finite space, or of the loop space
of a finite space are always locally finite. On the other side the
cohomology of the classifying space of a finite group whose order is
divisible by 2 does belong to any sub-category U_n.  One proves this
conjecture, modulo the additional hypothesis that all quotients of the
nilpotent filtration are finitely generated. This condition is used
when applying N. Kuhn's reduction of the problem. It is necessary to
do it to be allowed to apply Lannes' theorem on the cohomology of
mapping spaces.[N. Kuhn, On topologically realizing modules over the
Steenrod algebra, Ann. of Math. 141 (1995) 321-347].}

\maketitle

Un probl\`eme important en topologie est de savoir quand un module
instable sur l'alg\`ebre de Steenrod, ou une alg\`ebre instable,
est la cohomologie singuli\`ere d'un espace. L'exemple le plus
c\'el\`ebre est celui de la r\'ealisabilit\'e des alg\`ebres de
polyn\^omes comme cohomologie d'espaces. Ce probl\`eme, pos\'e par N.
E. Steenrod, est tr\`es li\'e au probl\`eme de l'invariant de Hopf
$1$ r\'esolu par J. F. Adams. Les techniques d\'evelopp\'ees au
cours des ann\'ees 80 et 90 par H. Miller et J. Lannes
ont d\'ebouch\'e sur les r\'esultats spectaculaires de B. Dwyer, H. Miller
 et
C. Wilkerson, entre autres, dans cette direction. Le pr\'esent
article \'etudie une autre instance de ce probl\`eme, qui a
\'et\'e soulev\'ee par N. Kuhn dans [3]. L'\'enonc\'e obtenu
prolonge un r\'esultat pr\'ec\'edent de l'auteur [8], r\'esultat
qui avait \'et\'e aussi conjectur\'e par Kuhn.

La cat\'egorie ${\cal U}$ des modules
instables sur l'alg\`ebre
de Steenrod admet une filtration
d\'ecroissante par des sous-cat\'egories pleines ${\cal N}il_k$, $k \geq 0$.
Ces cat\'egories sont d\'efinies comme suit.
La cat\' egorie $\Nil_k$ est la plus petite sous-cat\'egorie \'epaisse,
stable par limite
directe,  et qui contient $\Sigma^kM$ pour
tout
module instable $M$.
Chaque module instable admet donc
une filtration d\'ecroissante, que l'on appellera
filtration nilpotente du module. Si $M$ est un module instable
on notera $M_s$ son plus grand
sous-module instable $s$-nilpotent. Le quotient
$M_s/M_{s+1}$ est de la forme $\Sigma^sR_s(M)$, o\`u
$R_s(M)$ est un module instable r\' eduit, c'est-\`a-dire
ne contenant pas de suspension non-triviale.

 Cet article d\'emontre le th\'eor\`eme suivant~:

\begin{thm} Soit $X$ un espace, et soit $M$
la cohomologie singuli\`ere modulo $2$ de l'espace $X$.
Supposons que $M \in {\cal U}_n$ pour un certain
entier $n$, et que pour tout
$s$ le quotient
 $M_s/M_{s+1}$
ait un nombre fini de g\'en\'erateurs sur l'alg\`ebre
de Steenrod. Alors $M$ est localement
finie, c'est-\`a-dire limite directe
de modules finis sur l'alg\`ebre de Steenrod, {\it i.e.} $M \in {\cal U}_0$.
\end{thm}

En fait on d\'emontrera que chacun des quotients $M_s/M_{s+1}$ est
fini, et non nul seulement en degr\'e $s$ car il est r\'eduit. On en d\'eduit
que $M$ est localement fini car ~:

\begin{lem}
Un module instable $M$ tel que tous les quotients
$R_s(M)$ sont localement finis est localement fini.
\end{lem}

Ce lemme r\'esulte de l'exactitude de $T$, du fait que $T$ commute
aux suspensions, et de ce que $T(M)=M$ si $M$ est localement fini,
voir aussi [3].

Par comparaison, le cas de [8] est celui o\`u la filtration
nilpotente est finie, c'est-\`a-dire que le quotient $M_s/M_{s+1}$
est nul pour tout $s$ assez grand.

Une analyse pr\'ecise de la d\'emonstration
montre que l'on peut obtenir
des \'enonc\'es plus g\'en\'eraux.
A titre d'exemple, le suivant (que par pr\'ecaution
nous appelons conjecture)  semble \'a port\'ee
de main :

\begin{conj}
Soit $X$ un espace $2d$-connexe, $d \geq 1$. Soit $M$
sa cohomologie r\'eduite modulo $2$, supposons
que tous les cup-produits y soient triviaux. Supposons que
$M \in {\cal N}il_d$,   que $M/M_{2d} \in {\cal U}_1$,
 et que pour tout entier $s \leq d$ le
module instable $\bar T^s(M)$ soit de dimension finie
en chaque degr\'e, alors
$M \in {\cal N}il_{d+1}$.
\end{conj}

Il convient de noter que, relativement \`a la fitration de Krull, seul intervient le quotient
$M/M_{2d}$.

Soit $X$ un espace dont la cohomologie $M$ satisfait \`a
la condition de finitude
 impos\'ee ci-dessus sur les quotients de la filtration
nilpotente. Cette condition permet d'appliquer le th\'eor\`eme de
Lannes sur la cohomologie des espaces fonctionnels. En effet elle
implique que $T(M)$ est de dimension finie en chaque degr\'e. On
peut donc calculer la cohomologie de la cofibre $C(X)$ de $X \ra
{\rm map}(RP^\infty,X)$ et effectuer la r\'eduction de Kuhn [3].
C'est-\`a-dire que si $H^*X \in {\cal U}_n$, mais $H^*X \not \in
{\cal U}_{n-1}$, la cohomologie $ H^*C(X) \cong \bar T(H^*X) \in
{\cal U}_{n-1}$, mais $H^*C(X) \not \in {\cal U}_{n-2}$. On peut
donc raisonner par l'absurde et supposer qu'il existe un espace
$X$ dont la cohomologie est dans ${\cal U}_n$, mais non dans
${\cal U}_{n-1}$, et par it\'eration, se ramener au cas d'un
espace dont la cohomologie appartient \`a ${\cal U}_1$, mais pas
\`a ${\cal U}_0$ et montrer que ceci est impossible.

La conjecture, sans l'hypoth\`ese de finitude faite ci-dessus, ne
peut \^etre analys\'ee dans la cat\'egorie ${\cal U}$
et doit \^etre remplac\'ee par une conjecture portant
sur la structure de l'homologie comme module instable \`a droite.
Dans ce contexte on doit modifier la d\'efinition de la
 filtration de Krull dans le sens sugg\'er\'e par [4], {\it i.e.}
on doit consid\'erer une fitration d\'ecroissante sur
la cat\'egorie des modules instables \`a droite. Ceci pose aussi la question d'une extension du th\'eor\`eme de
 Lannes.

Une autre approche pour lever cette restriction de finitude
est d'appliquer les techniques de pro-espaces de F. Morel,
ceci a \'et\'e mis en oeuvre avec succ\`es par F-X. Dehon et G. Gaudens.

Il reste  la forme la plus g\'en\'erale de la conjecture de Kuhn
(\'evoqu\'ee plus haut). Nous reformulons ici l\'eg\`erement cette
conjecture, et pour ce faire le language des foncteurs [2], [7]
est plus commode. En fait un certain nombre d'exemples sugg\`erent
que la description  -et la technologie associ\'ee- qui suit
pourraient constituer une approche int\`eressante \`a cette
conjecture.

Rappelons
que l'on peut associer \`a un module instable r\'eduit un foncteur
analytique (limite
directe de foncteurs polyn\^omiaux) de la cat\'egorie des
espaces vectoriels de dimension finie sur le corps ${\bf F}_2$ dans
la cat\'egorie des
espaces vectoriels sur le corps ${\bf F}_2$, par exemple au module $F(1)$
est associ\'e le foncteur identit\'e
de degr\'e $1$. Ce foncteur
d\'etermine le module initial \`a localisation pr\`es.  On a en fait une \'equivalence
de cat\'egorie. Le module $R_s(H^*X)$
peut donc \^etre interpr\'et\'e
comme un tel foncteur ~:

\begin{conj} Soit $X$ un espace.
La suite des foncteurs $R_s(H^*X)$ est telle que~:
\begin{itemize}
\item soit il y en a au moins un foncteur dont
la s\'erie de Loewy (ou s\'erie des socles) est infinie et dont les facteurs de composition
ne sont pas born\'ees en degr\'e,
\item soit tous les foncteurs $R_s$ sont constants.
\end{itemize}
\end{conj}

En termes de modules instables cette conjecture 0.2 implique
que l'un au moins
des modules instables correspondants
a un nombre infini de g\'en\'erateurs.

Le fait que la s\'erie des socles d'un foncteur associ\'e \`a un
module est infinie est plus fort que le fait que le module  ne
soit pas de type fini ({\it i.e.} ai un nombre fini de
g\'en\'erateurs). La propri\'et\'e correspondante du module est
plus ennuyeuse \`a exprimer et moins naturelle. On renvoie \`a [9]
pour des d\'etails.

 A ce propos on notera que G. Gaudens a montr\'e
que la s\'erie des socles du foncteur
associ\'e \`a une alg\`ebre instable r\'eduite non concentr\'ee
en degr\'e z\'ero est infinie.

On peut pr\'eciser 0.2~:
Soit $d$ le plus petit entier
tel que $R_d$ soit
un foncteur non-constant et de degr\'e fini,  supposons $d > 1$. Alors on a~:

\begin{conj}
Un, au moins, des foncteurs $R_s$, $d<s<2d+1$
est de degr\'e non born\'e.
\end{conj}

La conjecture 0.1 ci-dessus est un argument en faveur de 0.2. Un
autre argument est l'exemple donn\'e par Kuhn dans [3] de la
filtration bar sur $K({\bf Z},3)$. Le second cran de cette
filtration fournit un espace $X$ tel que $R_1(H^*X)$ soit
isomorphe \`a $F(1)$, et donc \`a l'identit\'e comme foncteur,
mais $R_2(H^*X)$ a lui un nombre infini de g\'en\'erateurs et a,
en tant que foncteur, une s\'erie des socles infinie.

Le plan de la d\'emonstration du th\'eor\`eme 0.1 est le m\^eme
que dans [8]. Comme on l'a d\'ej\`a  dit, le cas de [8] correspond
\`a l'\'enonc\'e $0.1$  avec en plus l'hypoth\`ese que les
quotients $M_s/M_{s+1}$ sont nuls pour tout entier $s$ assez
grand.
 Rappelons tr\`es bri\`evement
l'id\'ee de [8]. On raisonne par l'absurde, et on suppose donc
qu'il existe un espace $X$ dont la cohomologie est dans ${\cal
U}_1$ mais n'est pas localement finie, {\it i.e.} n'appartient pas
\`a ${\cal U}_0$. On  montre alors que, dans la cohomologie d'un
certain espace de lacets it\'er\'es de $X$, la relation reliant le
cup-carr\'e aux op\'erations de Steenrod ne peut \^etre
satisfaite.

La diff\'erence fondamentale avec [8]  est que l'hypoth\`ese faite
dans cet article ne permet pas d'obtenir de zones d'annulation
dans la cohomologie des espaces de lacets associ\'es. On est
amen\'e
 \`a introduire
 des zones d'annulation modulo des termes
de degr\'e sup\'erieur de nilpotence.

Ceci fait qu'il est difficile d'obtenir un r\'esultat concernant
des complexes finis, dont on d\'eduirait celui cherch\'e, comme
cela est fait dans [8]. Ceci est n\'eanmoins possible, mais il y a
peu d'espoirs d'obtenir des \'enonc\'es g\'en\'eraux
satisfaisants.

H. J. Baues a sugg\'er\'e  la question suivante. Soit
$M$ un module instable, supposons
que $M \cong H^*X$. Quelles sont plus g\'en\'eralement les restrictions sur la structure de
$M$ impos\'ees par
les propri\'et\'es de la suite spectrale d'Eilenberg-Moore calculant
 la cohomologie de $\Omega X$, en particulier
impos\'ee par le fait que l'aboutissement de la suite spectrale
doit avoir une structure d'alg\`ebre instable?

Les sections $1$ et $2$ consistent d'abord en
 l'\'enonc\'e des r\'esultats concernant le comportement
des foncteurs $R_s$ sur la cohomologie
des espaces quand on passe de $X$ \`a $\Omega X$, puis
en des rappels sur la filtration nilpotente et la filtration de Krull. Dans la section 3
on d\'emontre ces r\'esultats
\`a l'aide de la suite spectrale d'Eilenberg-Moore. On
d\'emontre le th\'eor\`eme dans les sections $4$ et $5$.
La section $6$ donne un compl\'ement sur
la filtration nilpotente.

L'id\'ee essentielle de cet article a \'et\'e trouv\'ee lors d'un s\'ejour
au CRM \`a l'Univer\-sitat Aut\`onoma de Barcelona en Juin 1998. L'auteur tient \`a
remercier le groupe de topologie alg\'ebrique, et en particulier
J. Aguade et C. Broto, pour leur accueil chaleureux, Dagmar Meyer pour
ses remarques, et
le groupe de topologie alg\'ebrique de Tunis pour l'int\'er\^et
qu'il a manifest\'e
pour ce probl\`eme.

 {\bf On fera partout l'hypoth\`ese que
les modules instables consid\'er\'es sont de dimension finie en
chaque degr\'e.}

L'auteur tient \`a remercier le rapporteur pour lui avoir
signal\'e quelques ambiguit\'es, entre autres dans la d\'efinition
des classes $\alpha$, et dans la formulation d'un r\'esultat de
[3]. Il le remercie aussi pour lui avoir signal\'e diverses
r\'ef\'erences, en l'occurence celle de 0.2 dans [3] et celles de
propositions 2.4 et 2.5 dans la m\^eme r\'ef\'erence.

\section{La filtration nilpotente de la cat\'egorie ${\cal U}$, les foncteurs $R_s$
et les espaces de lacets}

On rappelle d'abord les deux d\'efinitions de la filtration
nilpotente sur la cat\' eg\-orie $\U$ des modules instables.
L'essentiel du mat\'eriel, concernant cette filtration, qui suit
ne pr\'etend pas \`a l'originalit\'e, il est soit explicitement
dans [6], [7], soit en est cons\'equence imm\'ediate (voir [3]).
Cependant les quelques ajouts (Propositions 1.8 \`a 1.12), faciles
eux aussi, sont n\'ecessaires \`a un traitement clair de la suite.
 Nous avons dans ce texte conserv\'e les
conventions de [7], on prendra garde au d\'ecalage dans la
d\'efinition des cat\'egories ${\cal N}il_k$ par rapport \`a [6].
Dans le contexte de [6] la cat\'egorie ${\cal U}$ est ${\cal
N}il_{-1}$. La notation choisie ici, celle de [7], pour ${\cal U}$
est ${\cal N}il_0$.

\begin{defn}
Soit un entier $s \geq 0$, la cat\' egorie $\Nil_s$ est la plus petite sous-cat\'egorie \'epaisse
stable par limite
directe, et qui contienne $\Sigma^sM$ pour
tout
module instable $M$. Un module qui appartient \`a $\Nil_s$ est dit
de degr\'e de nilpotence (au moins) $s$, ou (au moins)
$s$-nilpotent.
\end{defn}

On
remarquera qu'un module $s$-nilpotent est $(s-1)$-connexe.

La filtration nilpotente est d\'ecroissante et convergente.  La filtration de la cat\'egorie induit
sur chaque module instable une filtration d\'ecroissante, que l'on
appellera filtration nilpotente du module.

{\bf Si $M$ est un module
instable on notera $M_s$ son plus grand sous-module instable
$s$-nilpotent, cette notation sera conserv\'ee
dans tout l'article.}

\begin{defn}
 On dira qu'un \' el\' ement d'un module
instable
est de degr\' e
 de nilpotence au moins $s$, ou $s$-nilpotent, si le sous-module qu'il engendre est de degr\' e de nilpotence
au moins $s$. On dira qu'un \'el\'ement
est de degr\'e de nilpotence exactement $s$ s'il est de degr\'e
au moins $s$ et s'il n'est pas de degr\'e au moins $s+1$.
Un tel \'el\'ement
est n\'ecessairement non-nul.
\end{defn}

\begin{defn}
On note $\Sq_k$ l'op\' eration qui \`a un \' el\' ement $x$ de degr\' e $|x|$ d'un module instable $M$
associe la classe $\Sq^{|x|-k}x$ de degr\' e $2|x|-k$.
\end{defn}

Par convention  $\Sq_{k}x=0$
d\`es que $k>|x|$.
A cause de l'instabilit\'e on a $\Sq_{k}x=0$
d\`es que $k<0$.

Soit $M$ un module instable et soit $\Sigma^s(M)$.
 Soit $x \in M$, on note
$\sigma^sx$ la $s$-i\`eme suspension de $x$. On a~:
 $$\sigma^{s}(\Sq_k)( x))=(\Sq_{k+s})(\sigma^{s} x)
\  .$$

\begin{prop}{\rm[6]}\qua
Soit $M$ un module instable. Les conditions suivantes sont 
\'equivalentes~:
\begin{itemize}
\item $M$ est (au moins) $s$-nilpotent,
\item pour tout $x \in M$ et tout entier $k$ tel que
$0 \leq k < s$,
il existe un entier $c$,
d\' ependant de $x$ et $k$,
tel que $(\Sq_k)^c \, x=0$.
\end{itemize}
\end{prop}

Pour $s=0$ la condition est vide, on trouve donc bien ${\cal U}$.

C'est par la seconde condition de cette proposition que la
filtration nilpotente a \'et\'e introduite dans [6]. En fait, dans
les applications qui en ont \'et\'e donn\'ees jusqu'ici, c'est,
presque toujours, la premi\`ere condition ({\it i.e.} la
d\'efinition 1.1) qui a  \'et\'e utilis\'ee. Dans cet article  la
condition originelle de [6] joue un r\^ole important. Pour cette
raison, et parce que les indications de [6] sont un peu s\`eches,
on redonnera en fin de l'article la d\'emonstration du point cl\'e
de la proposition pr\'ec\'edente.

Soit $M$ un module instable
et $M_s$ son plus grand
sous-module instable $s$-nilpotent.  Chaque module instable admet
une filtration d\'ecroissante et convergente~:
$$
 M=M_{0} \supset M_{1} \supset  M_{2} \supset \ldots  M_{s} \supset
\ldots
$$
 Reprenant la notation de [3] on notera le quotient
$M_s/M_{s+1}$ sous la forme $\Sigma^sR_s(M)$, $R_s(M)$ est un
module instable {\bf r\' eduit}, c'est-\`a-dire ne contenant pas
de suspension non-triviale ([7] section 2.6). La d\'efinition de
$R_s$ est  naturelle et $R_s$ est un foncteur en $M$. On trouvera
dans  [3] une caract\'erisation intrins\`eque de la filtration
nilpotente.

Par ailleurs, comme
 les foncteurs $T$ et $\bar T$
de Lannes commutent aux suspensions et sont
exacts, ils respectent cette filtration. Par cons\'equent on a pour tout module instable $M$~:
$$
T(M_s)=T(M)_s \quad {\rm et} \quad \bar T(M_s)=\bar T (M)_s \  .
$$
Deux modules instables {\bf r\'eduits} $M$ et $N$ seront dits
fortement $F$-isomorphes si ils ont m\^eme ${\cal
N}il$-localisation ([7] section 6.3). On notera $M\cong_F N$. Ceci
revient \`a dire  qu'il existe un module instable $L$ tel que~:
\begin{itemize}
\item $L$ admet un monomorphisme $i$ dans $M$ et un monomorphisme $j$
dans $N$,
\item pour tout \'el\'ement non-nul $x$ dans $M$ il existe
un entier $c$ tel que $(\Sq_0)^cx $ est dans l'image de $i$ et non-nul,
pour tout \'el\'ement non-nul $x$ dans $N$ il existe
un entier $c$ tel que $(\Sq_0)^cx $ soit dans l'image de $j$ et non-nul.
\end{itemize}
Le module $L$ ci-dessus est n\'ecessairement r\'eduit.

Un monomorphisme $i~:~N \ra M$ est un $F$-isomorphisme fort si et seulement
si pour tout $x \in M$, $x \not = 0$, il existe un entier $c$ -d\'ependant
de $x$- tel que $\Sq_0^c(x) \in N$, et $\Sq_0^c(x) \not = 0$.

%\vskip 5mm

Les trois th\'eor\`emes
suivants d\'ecrivent le comportement de certains des foncteurs $R_s$
quand on passe d'un espace $X$ \`a son espace de lacets $\Omega X$. Ils seront
d\'emontr\'es dans la section 3, ce sont eux qui permettent de d\'emontrer la conjecture.
Apr\`es leur \'enonc\'e, les sections 1 et 2 sont consacr\'es aux pr\'eliminaires alg\'ebriques
\`a leur d\'emonstration. Ainsi qu'il a \'et\'e dit dans l'introduction
ces pr\'eliminaires se trouvent pour une part substantielle dans les r\'ef\'erences,
n\'ean\-moins un certain nombre de pr\'ecisions apport\'ees ici
seraient totalement hors contexte, pour le non-expert, si
des rappels n'\'etaient pas effectu\'es.

\begin{thm}
Soit $X$ un espace $1$-connexe tel que $H^*X$ est de dimension
finie en chaque degr\'e
et que $H^*X \in {\cal N}il_d$, $d > 1$. Alors $H^*\Omega X \in {\cal N}il_{d-1}$
et
pour $d-1 \leq s < 2d-2$ on a~:
$$
R_s(H^*X) \cong_F R_{s-1}(H^*\Omega X) \  .
$$
\end{thm}

On peut se demander si il n'y a pas, en fait, isomorphisme.

\begin{thm}
Soit $X$ un espace $1$-connexe tel que $H^*X$ soit de dimension
finie en chaque degr\'e
et que $H^*X \in {\cal N}il_d$, $d > 1$. Alors $H^*\Omega X \in {\cal N}il_{d-1}$
et
le module
instable $R_{2d-2}(H^*\Omega X)$ est fortement $F$-isomorphe \`a un module
instable $E$
donn\'e par une extension
de la forme~:
$$
 \{0\} \ra R_{2d-1}(H^*X) \ra E \ra M \ra \{0\} \ ,
 $$
o\`u $M$ est un sous-module de $ (R_{d}(H^*X))^{\otimes 2}$.

 \end{thm}

Soit $F_{-2}$ le deuxi\`eme terme de
la filtration d'Eilenberg-Moore de $H^*\Omega X$ (voir section 2).
Si $d=1$ les \'enonc\'es pr\'ec\'edents sont remplac\'es par~:

\begin{thm}
Soit $X$ un espace $1$-connexe tel que $H^*X$ est de dimension
finie en chaque degr\'e
et que $H^*X \in {\cal N}il_1$. Alors le module
instable $R_{0}(F_{-2}(H^*\Omega X))$ est fortement $F$-isomorphe \`a un module $E$
donn\'e par
une extension de la forme~:
$$
 \{0\} \ra R_{1}(H^*X) \ra E \ra M\ra \{0\} \ ,
 $$
o\`u $M$ est un sous-module de $ (R_1(H^*X))^{\otimes 2} $.
 \end{thm}

Les propositions 1.8 et 1.9 d\'ecrivent le comportement
des foncteurs $R_s$ par rapport aux suites exactes.

\begin{prop}
Soit $M$ un module instable, et soit $d$ un entier
donn\'e, $k \geq 1$. On suppose qu'il existe
une suite exacte~:
$$
0 \ra K \ra M \ra N \ra 0 \ ,
$$
avec $K \in \Nil_{k}$.
Alors si $s < k$ on a~:
$$R_s(M) \cong R_s(N)
\  .
$$
\end{prop}

Il suffit de d\'emontrer que l'application
$M/M_k \ra N/N_k$ est bijective. Elle est surjective
par construction, il suffit de donc montrer qu'elle est injective.
 La difficult\'e vient de ce que
le foncteur $M \mapsto M_k$ n'est pas exact \`a droite. Mais un \'el\'ement
$x \in M/M_k$ d'image nulle dans $N/N_k$ se rel\`eve en un \'el\'ement
$y \in M$ dont l'image dans $N$ appartient au sous-module
$N_k$. Il faut montrer que $y \in M_k$. Pour ce faire il suffit de montrer
que
$(\Sq_i)^cy$ est nul d\`es que $c$ est assez grand, pour $0 \leq i <k$.
Soit $z$ l'image de $y$ dans $N$, $z \in N_k$.

On a donc~:
$$
(\Sq_i)^cz
=0, \qquad 0 \leq i <k, \qquad {\rm si} \quad c \geq c_0  \ .
$$
Chacun des \'el\'ements $(\Sq_i)^{c_0}y$, $0 \leq i <k$ appartient donc
\`a $C$ qui est $k$-nilpotent. Le r\'esultat suit.

\begin{prop}
Soit $M$ un module instable, et soit $d$ un entier
donn\'e, $d \geq 1$. On suppose qu'il existe
une suite exacte~:
$$
0 \ra N \ra M \ra C \ra 0 \ ,
$$
avec
\begin{itemize}
\item $N \in \Nil_d$,
\item $C \in \Nil_{2d}$.
\end{itemize}
Alors~:
\begin{itemize}
\item le module $M$ est dans $\Nil_d$,
\item  si $d \leq s <2d$ on a~:
$$R_s(N) \cong_F R_s(M) \  ,
$$
\item  $R_{2d}(M) \cong_F E$  o\`u $E$ est
donn\'e par
 une extension de la forme~:
$$
\{0\} \ra R_{2d}(N) \ra E  \ra L \ra \{0\}   \  ,
$$
 $L$ d\'esignant  un sous-module de $R_{2d}(C)$.
 \end{itemize}
\end{prop}

{\bf D\'emonstration~:}

Une application de modules instables $f~:~N \ra M$ respecte la filtration nilpotente
et induit donc pour tout $s$ une application $f_s~:~N_s \ra M_s$
et une application $\bar f_s~:~R_s(N) \ra R_s(N)$. Le foncteur
$M \mapsto M_s$ est exact \`a gauche.  Il est facile de voir
que si $f$ est un monomorphisme
il en est de m\^eme pour chaque $\bar f_s$.

Le premier \'enonc\'e de la proposition
est clair.

%\vskip 5mm

Passons \`a  la seconde partie, soit $s < 2d$.

Notons $i$ l'injection de $N$ dans $M$ et $p$
la projection de $M$ sur $C$.
Consid\'erons le diagramme o\`u les notations sont \'evidentes~:
$$\diagram{
 \{0\} \ra & N_{s+1}  &\buildrel i
\over{\ra}  & M_{s+1} &\buildrel  p \over{\ra}  & C_{} \cr &
\vfl{}&& \vfl{} && \vfl{}  \cr \{0\} \ra & N_{s}  &\buildrel i
\over{\ra}  & M_{s}  &\buildrel  p \over{\ra}  & C_{} \cr &
\vfl{}&& \vfl{} && \vfl{}  \cr \{0\} {\ra}  & \Sigma^{s}R_{s}(N)
&\buildrel \bar i \over{\ra}  & \Sigma^{s}R_{s}(M) &\buildrel \bar
p \over{\ra} & \{0\} \cr }
$$
On notera que $C_{s+1}=C_s=C$.

 Par hypoth\`ese le conoyau de l'injection
$i_s~:~N_s \ra \ M_s$ est dans ${\cal N}il_{2d}$ et $C_{s+1}=C_s=C$. Soit $\bar x$
une classe non-nulle dans $M_s/M_{s+1}$ et soit
$x \in M_s$, $x \not = 0$ un rel\`evement
dans $M_s$. Soit  $y$ l'image de $x$ dans $C$. Comme $s < 2d$ il existe un
 entier
 $c$ -d\'ependant de $x$- tel que
$(\Sq_s)^cy=0$. Donc la classe $(\Sq_s)^cx$ est d'image nulle dans
$C$ et si $k \geq c$   $(\Sq_s)^kx \in i(N)$.

Soit $\sigma^{-s}\bar x \in R_s(M)$
 la $s$-i\`eme
d\'esuspension de l'image de
$$\bar x \in M_s/M_{s+1} \cong \Sigma^s R_s(M) \  .
$$
On a~:
$$\sigma^{-s}(\Sq_s)^k(\bar x))=(\Sq_0)^k(\sigma^{-s}\bar x) \in (\bar(R_s(N))
\  ,$$
si $k \geq c$.
On a donc d\'emontr\'e
que $(\Sq_0)^k(\sigma^{-s}\bar x)$ est non-nulle et dans l'image de $\bar i_s$
pour tout $k \geq c$. Ceci
donne la seconde partie de la proposition.

%\vskip 5mm

 Passons \`a la troisi\`eme partie, et faisons donc $s=2d$.
La d\'emonstration
va occuper la fin de cette section. Si on a une suite exacte~:
$$
\{0\} \ra N \ra M \ra C \ra \{0\}
$$
le complexe~:
$$
\{0\} \ra R_{2d}(N) \ra R_{2d}(M) \ra R_{2d}(C) \ra \{0\}
$$
n'est pas exact au centre et \`a droite. En fait le complexe de
foncteurs [2]~:
$$
\{0\} \ra f(R_{2d}(N)) \ra f(R_{2d}(M)) \ra f(R_{2d}(C))
$$
est exact.
 Ceci se traduit, en termes de modules instables,
 par la proposition suivante~:

\begin{prop} Il
existe un diagramme commutatif de modules instables o\`u la ligne
sup\'erieure
est exacte~:
$$\diagram{
\{0\} \ra & R_{2d}(N) \ra & K \ra & L \ra & \{0\} \cr
{} & \vfl{{\rm Id}}& \vfl{i} & \vfl{j} & {}\cr
\{0\} \ra & R_{2d}(N) \ra & R_{2d}(M) \ra & R_{2d}(C) \cr
}$$
et o\`u $i$ est un $F$-isomorphisme fort, et $j$ est
un monomorphisme.
\end{prop}

D\'emontrons cette propri\'et\'e.
Comme le foncteur $M \mapsto M_{2d}$ est exact \`a gauche
on a \'evidemment un complexe exact \`a gauche et au centre~:
$$
0 \ra N_{2d} \ra M_{2d} \ra C_{2d} \ ,
$$
et donc un complexe~:
$$
0 \ra
R_{2d}(N) \ra R_{2d}(M) \ra R_{2d}(C)  \  ,
$$
o\`u, {\it a priori}, on sait seulement que la premi\`ere application est injective. Avec
les m\^emes notations que plus haut, on a un
diagramme commutatif~:
$$\diagram{
 \{0\} \ra & N_{2d+1}  &\buildrel i
\over{\ra}  & M_{2d+1} &\buildrel  p \over{\ra}  & C_{2d+1} \cr &
\vfl{}&& \vfl{} && \vfl{}  \cr \{0\} \ra & N_{2d}  &\buildrel i
\over{\ra}  & M_{2d}  &\buildrel  p \over{\ra}  & C_{2d} \cr &
\vfl{}&& \vfl{} && \vfl{}  \cr \{0\} {\ra}  & \Sigma^{2d}R_{2d}(N)
&\buildrel \bar i \over{\ra}  & \Sigma^{2d}R_{2d}(M) &\buildrel
\bar p \over{\ra} & \Sigma^{2d}R_{2d}(C) \cr }
$$
La ligne inf\'erieure est exacte au milieu
au sens  donn\'e par le lemme suivant~:

\begin{lem} Si un \'el\'ement $x \in \Sigma^{2d}R_{2d}(M)$ est dans
le noyau de $\bar p$, alors pour
tout $k$
assez grand $(\Sq_{2d})^kx$ est dans l'image de $\bar i$.
\end{lem}

La d\'emonstration
est identique \`a celle donn\'ee plus haut et laiss\'ee au lecteur.

Pour d\'emontrer la proposition
il  reste \`a construire le diagramme~:
$$\diagram{
\{0\} \ra & R_{2d}(N) \ra & K \ra & L \ra & \{0\} \cr
{} & \vfl{Id}& \vfl{i} & \vfl{j} & {}\cr
\{0\} \ra & R_{2d}(N) \ra & R_{2d}(M) \ra & R_{2d}(C) \cr
}$$
avec les propri\'et\'es requises.

Le lemme suivant permet de construire $K$~:

\begin{lem} Soit $H$ un module instable r\'eduit, et soit
$J$ un sous-module. Alors il existe un sous-module $H'$
de $H$ tel que~:
\begin{itemize}
\item $H' \cong_F H$,
\item si $x \in H'$ et $\Sq_0x \in J$ alors $x \in J$.
\end{itemize}
\end{lem}

La d\'emonstration proc\`ede en deux \'etapes. D'abord si $H$
a un nombre fini de g\'en\'erateurs sur
l'alg\`ebre de Steenrod on montre que
$(\Sq_0)^k(H)+J$ convient d\`es que $k$ est assez grand. En
effet, soit $J_h$ le sous-module
de $H$ constitu\'e
par les \'el\'ements $x \in H$ tels que
$(\Sq_0)^hx \in J$, on
a $J_h \subset J_{h+1} \subset \ldots$. Comme $J$ a
un nombre fini de g\'en\'erateurs
il existe un entier $k$ tel que $J_k=J_{k+1}=\ldots$. On
v\'erifie que $(\Sq_0)^k(H) \cap J_k \subset J$ et
que  $(\Sq_0)^k(H)+J$ convient (on utilise l'injectivit\'e
de $\Sq_0$). Puis un argument de limite directe convenable, utilisant
le fait que les modules consid\'er\'es sont de dimension finie
en chaque degr\'e permet de conclure. En fait
on ne consid\'erera que des modules ayant
un nombre fini de g\'en\'erateurs
et on n'aura pas besoin de cette g\'en\'eralisation.

Pour terminer la d\'emonstration de
la proposition on applique le lemme avec $H=R_{2d}(M)$ et $J=R_{2d}(N)$
Le module $L$ est alors  le quotient de $K$
par l'image de $R_{2d}(N) $.

En fait le lemme 1.12 peut \^etre pr\'ecis\'e en utilisant les
notations de [2], et notamment les foncteurs $m$ et $f$. On montre
que l'on a un diagramme commutatif~:
$$\diagram{
\{0\}& \ra & R_{2d}(N) &\ra & R_{2d}(M)& \ra & L &\!\!\!\!\!\ra & \{0\} \cr
{} & &\vfl{}& &\vfl{i} && \vfl{j} & &{}\cr
\{0\}& \ra & m \circ f(R_{2d}(N))& \ra & m \circ f(R_{2d}(M)) &\ra & m \circ f(R_{2d}(C))& \cr
}$$
o\`u chaque ligne est exacte et $j$ a un noyau dans ${\cal N}il_1$.

%\vskip 1cm

\section {La filtration de Krull de la cat\'egorie {\cal U}}
Rappelons la d\'efinition de la filtration de Gabriel-Krull sur la
cat\'egorie ${\cal U}$ (voir [1]).

Un module instable est localement fini s'il est limite directe de
modules instables finis. La sous-cat\'egorie pleine, ${\cal B}$ des modules
localement finis est le premier terme de la filtration de
Gabriel-Krull sur la cat\'egorie ${\cal U}$. Rappelons la
d\'efinition de cette filtration, ${\cal U}_{0}$ est la plus petite
classe de Serre de ${\cal U}$, stable par limite directe et
contenant tous les objets simples de ${\cal U}$.
 Comme ceux-ci sont de la forme
$\Sigma^n{\bf F}_{2}$, ${\cal U}_0$ s'identifie
 \`a ${\cal B}$. Cette construction
s'\'etend
\`a toute cat\'egorie ab\'elienne ${\cal A}$. Supposons que ${
\cal  U}_{n}$ ait \'et\'e d\'efinie, d\'efinissons alors ${\cal  U}_{n+1}$
comme \'etant la sous-cat\'egorie pleine de ${\cal U}$ d\'efinie
comme suit. Dans la cat\'egorie quotient  ${\cal  U}/{
\cal  U}_{n}$ on consid\`ere la plus petite sous-cat\'egorie
\'epaisse qui est stable par limite directe et qui contient tous les objets
simples de ${\cal U}/{\cal U}_{n}$. Alors, un objet $M$ de ${\cal
U}$ appartient \`a ${\cal U}_{n}$, si et seulement si en tant
qu'objet de la cat\'egorie ab\'elienne ${\cal U}/{\cal U}_{n}$, il
appartient
\`a la sous-cat\'egorie  $({\cal  U}/{\cal  U}_{n})_0$.

\begin{thm} {\rm([7] section 6.2)}\qua La plus petite sous-cat\'egorie de
${\cal U}$ contenant toutes les cat\'egories ${\cal U}_n$,
et stable par limite directe, est ${\cal U}$ elle m\^eme.
\end{thm}

Voici une caract\'erisation de la filtration de Krull
\`a l'aide
de $\bar T$.

\begin{thm} {\rm([7] 6.2.4)}\qua Soit $M$ un module
instable. Alors
  $M$ appartient \`a ${\cal  U}_{n}$ si et seulement si
 $\bar T^{n+1}M$ est trivial.
 \end{thm}

On  d\'ecrit maintenant en d\'etails les objets de ${\cal U}_1$,
ce r\'esultat est dans [8], avec l'hypoth\`ese suppl\'ementaire
que le module a un nombre fini de g\'en\'erateurs. La
d\'emonstration, ne change en rien, et est donn\'ee rapidement.

\begin{prop}Soit un module instable $M$ qui appartient
\`a ${\cal U}_1$, mais qui n'appartient pas \`a
${\cal U}_0$.
Soit
$K=\bar T(M)$, alors $K$ est localement fini et non-trivial,
et
il existe une suite
exacte~:
$$
0 \ra L \ra M \ra K \otimes F(1) \ra  L' \ra 0 \  ,
$$
o\`u $L$ et $L'$ sont localement finis.
\end{prop}

Par d\'efinition de $\bar T$ on a une application
$$
M \ra K \otimes \bar H^* RP^\infty.
$$
Le noyau de cette application est localement fini, notons le $L$.
En effet, la d\'efinition de $\bar T$ donne, par adjonction, $\bar
T L = \{0\}$. On a vu plus haut que cette condition caract\'erise
les modules localement finis.
 On d\'emontre
 comme dans [8],
 que cette application est \`a valeurs dans
$K \otimes F(1)$.

Pour montrer  que le conoyau $L'$ est localement fini, affirmation
la plus facile, il suffit d'appliquer le foncteur $\bar T$
\`a la suite exacte~:
$$
0 \ra L \ra M \ra K \otimes F(1) \ra L' \ra 0 \ ,
$$
ce qui donne $\bar T(L')={0}$.

On va maintenant d\'ecrire la filtration nilpotente sur
un tel module. Le r\'esultat suivant est une cons\'equence
\'el\'ementaire de la d\'efinition de
cette filtration.

\begin{prop} Soit $M$ un module instable localement fini. Le
sous-module
$M_s$ de $M$ est le sous-module des \'el\'ements de degr\'e sup\'erieur
ou \'egal \`a $s$.
\end{prop}

La d\'emonstration est laiss\'ee en exercice au lecteur (voir
aussi  [3]). Dans ce cas  la filtration nilpotente s'identifie
donc  \`a la filtration par le degr\'e.

Si $M \in {\cal U}_1$ on applique $1.14$. Si le module $K$ est
localement fini c'est encore un exercice facile de v\'erifier que
la filtration nilpotente sur $K \otimes F(1)$ est induite, par
produit tensoriel par $F(1)$, par la filtration par le degr\'e sur
$K$. Ceci est d\'emontr\'e en d\'etails et en plus grande
g\'en\'eralit\'e dans [3], et est d'ailleurs une cons\'equence
imm\'ediate de [6]. On en d\'eduit donc que~:
\begin{cor}
Si $K$ est localement fini on a~:
$$
R_s(K \otimes F(1)) \cong  K_s \otimes F(1) \  ,
$$
$K_s$ \'etant compris comme un module instable concentr\'e en degr\'e z\'ero.
\end{cor}

Ceci d\'ecrit, par restriction, la
filtration  sur le quotient $M/L$.

Le r\'esultat
suivant est corollaire des deux propositions pr\'ec\'edentes.

\begin{cor} Soit $M$ un module instable qui appartient \`a  ${\cal U}_1$.
Supposons que  $M$ doit $k$-connexe.
Soit alors $s$ un entier tel que $ s \leq k$, le module $R_s(M)$ est trivial
dans les degr\'es qui ne sont pas une puissance de $2$.
\end{cor}

La condition de connexit\'e est l\`a pour garantir
que $L$ est $k$-nilpotent, et que
l'on peut appliquer $1.8$.

On aura aussi consid\'erer des modules r\'eduits
qui sont produit tensoriel de deux modules dans ${\cal U}_1$
et des extensions de tels modules
par des modules dans ${\cal U}_1$. La proposition suivante est
 vraie, par construction, pour ces modules.
Ces modules sont dans la cat\'egorie ${\cal U}_2$  [FS], [7], et
 il est int\`eressant
de  donner la propostion dans son contexte g\'en\'eral, ne serait ce que parce
que pour montrer
que le r\'esultat n'est pas fortuit.

\begin{prop} {\rm[FS]}\qua Soit $M$ un module instable r\'eduit qui appartient \`a  ${\cal U}_2$.
Supposons que $M$ soit $0$-connexe.
Alors $M$ est trivial
dans les degr\'es qui ne sont pas de la forme $2^h$, $h \geq 0$, ou $2^i+2^j$,
$i>j\geq 0$.
\end{prop}

La condition de connexit\'e est l\`a pour ne pas
ne pas  avoir \`a inclure $0$ dans la liste des degr\'es.

Les \'enonc\'es qui suivent sp\'ecialisent la proposition 1.9
au cas consid\'er\'e dans l'article.

\begin{prop}
Soit $M$ un module instable, et soit $d$ un entier
donn\'e, $d \geq 1$. On suppose qu'il existe
une suite exacte~:
$$
0 \ra N \ra M \ra D \ra 0 \ ,
$$
avec
\begin{itemize}
\item $N \in \Nil_d$,
\item $D \in \Nil_{2d}$,
\item la connectivit\'e de $M$ est sup\' erieure ou \'egale \`a $2d$,
\item le quotient  $N/N_{2d}$,
 est dans ${\cal U}_1$,
\item le module  $R_{2d}(N)$ est dans ${\cal U}_1$,
le module  $R_{2d}(D)$ est dans ${\cal U}_2$.

\end{itemize}
Alors~:
\begin{itemize}
\item le module $M$ est dans $\Nil_d$ et donc  $R_s(M)=\{0\}$ si $s <d$.
,
\item le module $R_k(M)$,
avec $d \leq k <
 2d$, est dans ${\cal U}_1$, connexe, et
donc trivial dans les degr\' es qui ne sont
pas  une puissance de $2$,
\item le quotient $R_{2d}(M)
$ est dans ${\cal U}_2$, connexe,  et donc
 trivial  dans les degr\'es qui ne sont
pas de la forme
$2^h$, $h \geq 0$ ou de la forme $2^h+2^j$, $h>j\geq 0$.
\end{itemize}
\end{prop}

%\vskip 2mm

{\bf D\'emonstration~:} La premi\`ere partie est claire. La seconde
r\'esulte de 1.9 et 2.5, utilisant
que le quotient  $N/N_{2d}$,
 est dans ${\cal U}_1$. La troisi\`eme r\'esulte de 1.9
et 2.6, utilisant
que le module  $R_{2d}(N)$ est dans ${\cal U}_1$ et connexe, et que
le module  $R_{2d}(D)$ est dans ${\cal U}_2$ et connexe.

Dans l'\'enonc\'e les conditions de connectivit\'e sont de
convenance. On pourrait les supprimer et remplacer les
cat\'egories ${\cal N}il_k$ par les cat\'egories $\bar {\cal
N}il_k$ [7], ou encore rajouter dans les conclusions, aux degr\'es
$2^h$, le degr\'e $0$, et aux degr\'es $2^j+2^h$, le degr\'e $0$.
Ceci entra\^inerait quelques modifications mineures dans la suite,
dans les \'enonc\'es et les d\'emonstrations.
 Il faudrait prendre garde \`a la convergence de la suite spectrale d'Eilenberg-Moore.
Ceci pourrait \^etre fait comme dans [7].

Un \'enonc\'e analogue a lieu pour $d=0$~:

\begin{prop}
Soit $M$ un module instable. On suppose qu'il existe
une suite exacte~:
$$
0 \ra N \ra M \ra C \ra 0 \ ,
$$
avec
\begin{itemize}
\item  $M$ est connexe,
\item le quotient  $N/N_{1}\cong R_0(N)$,
 est dans ${\cal U}_1$, le quotient
 $C/C_{1}\cong R_0(C)$,
 est dans ${\cal U}_2$.
\end{itemize}
Alors~:
\begin{itemize}
\item le quotient $M/M_{1}\cong R_0(M)$ est
 dans ${\cal U}_2$, connexe, et donc trivial dans les degr\'es qui ne sont
pas de la forme $2^h$, $ \geq 0$ ou de la forme
  $2^h+2^j$  $j > h \geq 0$.
\end{itemize}
\end{prop}

{\bf D\'emonstration~:} Comme plus haut.

\section {\hglue -5pt Applications \`a la suite
spectrale d'Eilenberg-Moore}

On va appliquer les r\' esultats pr\'ec\'edents
 dans le contexte de la suite spectrale d'Eilenberg-Moore.
Rappelons en les propri\'et\'es principales [Re], [10].
 Soit $X$ un espace $1$-connexe.

La cohomologie modulo
 $2$
de $\Omega X$   a une filtration naturelle d\' ecroissante et
convergente par des modules instables~:
$$
 \ldots \supset F_{-s} \supset F_{-s+1} \ldots \supset F_{-1} \supset F_{
0} \supset \{0\} \  .
$$
Le quotient  $\Sigma^d(F_{-d} / F_{-d+1})$
 est isomorphe, en tant qu'espace vectoriel gradu\'e
\`a la colonne $E_\infty^{-d,*}$  de la
suite spectrale d'Eilenberg-Moore.

 Le terme $E_2^{-d,*}$ de
la suite spectrale est  isomorphe
\`a ${\rm Tor}_{H^*X}^{-d}(\F,\F)$. Ce module est un sous-quotient
de $\bar H^*X^{\otimes d}$, et la formule de Cartan y
d\'etermine une structure de module instable.

Soit $r \geq 2$, les diff\'erentielles s'interpr\`etent comme des applications de
degr\'e z\'ero~: $$d_r~:~E_r^{s,*} \ra \Sigma^{r-1}(E_r^{s+r, *})
\  .
$$
Ce sont
 des applications de modules instables.
On a donc sur le terme $E_\infty^{-d,*}$ une seconde structure de
module instable. Elle est identique \`a la pr\'ec\'edente.

On va \'etudier le comportement des foncteurs $R_s$ quand on
passe de $X$ \`a $\Omega X$. On a~:

\begin{prop}
Soit $X$ un espace tel que $H^*X$ soit de dimension
finie en chaque degr\'e
et que $\bar H^*X \in {\cal N}il_d$, $d > 1$. Alors $ \bar H^*\Omega X \in {\cal N}il_{d-1}$
et
pour $d \leq s < 2d-1$ on a~:
$$
R_s(\bar H^*X) \cong_F R_{s-1}(\bar H^*\Omega X) \  .
$$
\end{prop}

On a \'evidemment $R_s(\bar H^*X) =\{0\}$ si $s < d$, et $R_s(\bar H^*\Omega
X) =\{0\}$ si $s < d-1$.

%\vskip 2mm

{\bf D\'emonstration~:}
On va comparer $ H^*X$, $F_{-1}$ et $H^* \Omega X$. Pour ce faire on commence
par consid\'erer l'homomorphisme de coin, dont
on rappelle qu'il est induit par
l'application canonique d'\'evaluation $$\Sigma \Omega X \ra X
\  .
$$
Il se d\'ecompose comme suit~:
$$
\bar H^*X \ra QH^*X \ra \Sigma F_{-1} \ra \Sigma H^*\Omega X \  ,
$$
o\`u $QH^*X \cong \bar H^*X /(\bar H^*X)^2$.

Commen\c cons par examiner le conoyau de $F_{-1} \ra H^*\Omega X$.
Le module instable $\Sigma^s (F_{-s}/F_{-s+1})$ est un
sous-quotient de ${\rm Tor}_{H^*X}^{-s}(\F,\F)$, qui est lui
m\^eme un sous-quotient de $\bar H^*X ^{\otimes s}$. Or, il
r\'esulte de [6] que
 le produit tensoriel d'un module $L$ dans $\Nil_u$ par un
module $L'$ dans $\Nil_v$ est dans $\Nil_{u+v}$. On en d\' eduit
que $${\rm Tor}_{H^*X}^{-s}(\F,\F)$$ est dans $\Nil_{sd}$.
Le
sous-quotient it\'er\'e $E_\infty^{-s,*}$ est dans $\Nil_{sd}$.
Donc le quotient $F_{-s}/F_{-s+1}$ est dans $\Nil_{sd-s}$.
Le quotient $ H^*\Omega X/F_{-1}$ admet une filtration d\'
ecroissante dont les termes sont les modules $F_{-i}/F_{-2}$, $i \geq 2$. Changeant $-i$ en
$i$, $i \geq 2$ on obtient une filtration croissante dont le $i$-i\`eme quotient
 est dans $\Nil_{i(d-1)}$. Comme les cat\' egories
$\Nil_k$ sont des classes de Serre,
 et sont stables par limite directe, on en conclut
que le quotient $ H^*\Omega X/F_{-1}$ est $(2d-2)$-nilpotent.

On peut  appliquer la proposition $1.9$
\`a~:
$$
\{0\} \ra  F_{-1} \ra \bar H^*\Omega X \ra \bar H^*\Omega X/F_{-1} \  ,
$$
et en conclure que~:
$$
R_s( F_{-1}) \cong_F R_s( \bar H^*\Omega X) \ 
$$
pour $s <2d-2$.

Le noyau de $\bar H^*X \ra QH^*X$ est \'egal \`a
$(\bar H^*X)^2$ et est $2d$-nilpotent.
Il en r\'esulte (Proposition 1.8) que l'application~:
$$
R_s(\bar H^*X) \ra R_s(QH^*X)
$$
est un isomorphisme pour $s<2d$. Le noyau
de $QH^*X \ra \Sigma F_{-1}$ est l'image des diff\'erentielles
dans $E_2^{-1,*} \cong  QH^*X$. La description
des diff\'erentielles donn\'ee plus haut
et le fait que la colonne $E_r^{-s,*}$ soit
$sd$-nilpotente implique que cette image est
$2d$-nilpotente, $d \geq 1$.
Il en r\'esulte (Proposition 1.8) que l'application~:
$$
R_s(H^*X) \ra R_s(
\Sigma
F_{-1})
$$
est un isomorphisme pour $s<2d$.

L'ensemble de ces r\'esultats donne la proposition.

%\vskip 5mm

On \'etudie maintenant le cas de $R_{2d-
2}(H^*\Omega X)$. Le r\'esultat
suivant est cons\'equ\-ence de la proposition $1.9$~:

\begin{prop}
Soit $X$ un espace. On suppose que
$H^*X$ est $d$-nilpotent,  $d > 1$. Alors~:
\begin{itemize}
\item  $R_{2d-2}(\bar H^* \Omega X)$ est fortement $F$-isomorphe \`a un module $E$ donn\'e par
 une extension de la forme~:
$$
\{0\} \ra R_{2d-1}(\bar H^*X) \ra E  \ra L \ra \{0\}   \  ,
$$
o\`u $L$ est un sous-module de $R_{d}(\bar H^*X)^{\otimes 2}$
 \end{itemize}
Si $d=1$ on a un r\'esultat analogue en substituant
$F_{-2}$ \`a $H^* \Omega X$.
\end{prop}

{\bf D\'emonstration~:}
On commence par observer que
$$  R_{2d-2}(F_{-2}) \cong_F R_{2d-2}(\bar H^* \Omega X) \  .
$$
Ceci r\'esulte de $1.9$ et du fait que $H^* \Omega X / F_{-2}$ est $3d-3$-nilpotent, d\'emontr\'e comme plus haut.
Il reste donc \`a analyser $R_{2d-2}(F_{-2})$,
ce qu'on fait par  d\'evissage et application de $1.9$
en consid\'erant la suite exacte~:
$$
\{0\} \ra F_{-1} \ra F_{-2} \ra F_{-2}/F_{-1} \ra \{0\} \  .
$$
Le module instable  $\Sigma F_{-1}$ est -\`a suspension
pr\`es- l'image de l'homomorphisme de coin (voir la d\'emonstration
qui pr\'ec\`ede).
Le module instable $\Sigma^2 (F_{-2}/F_{-1})$ est un sous-quotient
de ${\rm Tor}_{H^*X}^{-2}(\F,\F)$, qui est lui m\^eme un
sous-quotient de $\bar H^*X^{\otimes 2}$.
Du fait que les diff\'erentielles ont des images
au moins $3d$-nilpotentes
et de la proposition
$1.8$ et on d\'eduit
que $R_{2d}(\Sigma^{2}(F_{-2}/F_{-1}))$ est isomorphe \`a
un sous-module de $R_d(\bar H^*X)^{\otimes 2}$.

Le r\'esultat suit.

%\vskip 5mm

Les \'enonc\'es suivants adaptent et pr\'ecisent
 les \'enonc\'es pr\'ec\'edents au
contexte, ils en sont cons\'equence directe.

\begin{thm} Soit $d$ un entier strictement sup\'erieur
$1$, soit  $X$ un espace, et soit
$\bar H^*X$ sa cohomologie r\'eduite. Supposons  que~:
\begin{itemize}
\item $\bar H^*X$ appartient \`a ${\cal N}il_d$,
\item la connectivit\'e de $\bar H^*X$ est strictement sup\'erieure \`a $2d$,
\item le quotient $H^*X/(H^*X)_{2d}$ est dans ${\cal U}_1$.
\end{itemize}
Alors
\begin{itemize}
\item $\bar H^*\Omega X$ appartient \`a ${\cal N}il_{d-1}$,
\item la connectivit\'e de $\bar H^*\Omega X$ est
sup\' erieure \`a $2d-2$,

\item le module $R_s(\bar H^* \Omega X)$, $d-1 \leq s <
 2d-2$, est dans ${\cal U}_1$, connexe,  et donc trivial  dans les degr\'es qui ne sont pas de la forme
$2^h$, $h \geq 0$,
\item le module
$R_{2d-2}(\bar H^*\Omega X)$
 est dans ${\cal U}_2$, connexe,  et donc trivial dans les degr\' es qui ne sont
pas de la forme $2^h$, $h \geq 0$, ou $2^h+2^j$, $h>j \geq 0$.
\end{itemize}
\end{thm}

La d\'emonstration est cons\'equence de 3.1, 3.2, 2.5, et 2.6.

%\vskip 5mm
Un \'enonc\'e particulier est n\'ecessaire dans le cas,
exclu
ci-dessus, o\`u $d=1$. L'\'enonc\'e est similaire, mais le r\'esultat
ne s'applique qu'au sous-module $F_{-2}$ de
$\bar H^*\Omega X$~:

\begin{thm} Soit $X$ un espace, et soit
$\bar H^*X$ sa cohomologie r\'eduite. Supposons que~:
\begin{itemize}
\item $\bar H^*X$ appartient \`a ${\cal N}il_1$,
\item la connectivit\'e de $\bar H^*X$ est sup\'erieure \`a $2$,
\item le module $R_1(\bar H^*X)$
est dans ${\cal U}_1$.
\end{itemize}
Alors
\begin{itemize}
\item le quotient $
R_0(F_{-2}(\bar H^* \Omega X))
$ est dans ${\cal U}_2$, connexe, et donc trivial
dans les degr\'es qui ne sont pas de la forme
$2^h$, $h \geq 0$ et $2^h+2^j$, $j >h \geq 0$.
\end{itemize}
\end{thm}

La d\'emonstration est identique \`a celle faite plus haut
mais on doit se restreindre \`a $F_{-2}$.

\section         {Construction de classes dans la cohomologie
des espaces de lacets }

On va maintenant donner le d\'ebut de la d\'emonstration du
th\'eor\`eme $0.1$. La r\'eduction de Kuhn  permet de supposer
 qu'il existe un espace
$X$ dont la cohomologie appartienne \`a ${\cal U}_1$,
mais pas \`a ${\cal U}_0$. Cette r\'eduction s'effectue comme
dans le cas pr\'ec\'edent, le seul probl\`eme est d'assurer les
conditions d'application du th\'eor\`eme de Lannes
sur la cohomologie des espaces fonctionnels.

On peut appliquer le th\'eor\`eme de Lannes pour la raison
suivante. On a suppos\'e que chaque quotient de la filtration
nilpotente a un nombre fini de g\'en\'erateurs, ceci implique que
$H^*X$ est de dimension finie en chaque degr\'e. De plus
l'hypoth\`ese est pr\'eserv\'ee par application de $ T$. Ceci
implique que $T(H^*X)$ est de dimension finie en chaque degr\'e,
ce qui  permet de calculer la cohomologie de l'espace fonctionnel
\`a l'aide du th\'eor\`eme de Lannes.

La proposition 2.3 s'applique \`a la cohomologie de $X$ en posant
$M=\bar H^*X$ et $K=\bar T(\bar H^*X)$. Le module instable $K$ est non
trivial. Soit $d-1$ sa connexit\'e, {\bf l'entier $d$ d\'esignera
dor\'enavant, et pour
toute la suite cette constante}. Le module instable $K$ est
inchang\'e quand on remplace $H^*X$ par un sous-module $M'$ tel que
le quotient $H^*X/M'$ soit fini. Ceci implique, en particulier, que
$K$ est inchang\'e, quand on quotiente $X$ par un sous-complexe de
dimension finie. On peut donc supposer, par commodit\'e, la
connectivit\'e de $X$ aussi grande que l'on veut, donc on peut la
supposer sup\'erieure \`a $2d$.

On doit traiter d'abord le cas o\`u $d=0$. Dans ce cas, par
hypoth\`ese, et \`a cause de 2.5, le quotient de $H^*X$ par
son sous-module nilpotent maximal est non-nul seulement
 en degr\'e de la forme $2^h$. Le sous-module
nilpotent maximal est un id\'eal, le quotient
est donc aussi une alg\`ebre instable. Soit $x$ est un
\'el\'ement non-nul dans ce quotient. Toutes les puissances
$x^{2^h}$ sont non nulles. L'\'el\'ement $x^3$ est aussi
 non-nul, car $x^4$ est non-nul.
 Le degr\'e de $x^3$ n'est
pas une puissance de $2$, il y a donc une contradiction et
$d$ n'est pas nul.

Supposons donc $d \geq 1$, et introduisons, comme dans [8], des
classes dans la cohomologie de $X$ et de ses espaces de lacets
it\'er\'es. Rappelons que l'on identifie $F(1)$, comme \`a
l'ordinaire, avec le sous-module instable, de $H^*B{\bf Z}/2\cong
\F[u]$, engendr\'e par $u$ et qui admet pour base sur $\F$ les
\'el\'ements $u^{2^i}$. On peut appliquer la proposition 2.3 \`a
$\bar H^*X$. Consid\'erons alors une classe $\omega \in K^d$,
$\omega \not = 0$ et les classes $$ \omega  \otimes u^{2^{j}} \in
K \otimes F(1) \  .
$$
Faisons l'hypoth\`ese suivante. Soit ${\cal A}\omega$ le
module instable engendr\'e par $\omega$. On sait qu'il est fini
puisque c'est un sous-module de $K$ qui est localement fini.
Supposons ${\cal A}\omega$ nul en degr\'e sup\'erieur ou \'egal \`a
$h$. Soit un entier $k$ tel que $2^{k-1} \geq  h$,
 sous cette hypoth\`ese
la formule de Cartan montre que~:
$Sq^{2^{k+i}} (\omega \otimes u^{2^{k+i}})=\omega \otimes u^{2^{k+i+1}}$.

D'apr\`es 2.3 il existe un entier positif $k_0$, d\'ependant de $\omega$,
tel que
$ \omega  \otimes u^{2^{j}} \in K \otimes F(1)$
se rel\`eve \`a  $\bar H^*  X$ d\`es que $j \geq k_0$.
On notera $\alpha_{i,d}$
ces rel\`evements, ils sont de degr\'e $2^{k_0+i}+d$. {\it A priori} ils ne sont d\'efinis
que modulo un \'el\'ement localement fini, mais on peut
les d\'eterminer de mani\`ere univoque  en choisissant $\alpha_{0,d}$ et en
supposant qu'ils v\'erifient
$Sq^{2^{k_0+i}} (\alpha_{i,d})=\alpha_{i+1,d}$.

{\bf L'entier $k_0$ restera  fixe dans toute la suite de la
d\'emonstration.}

Quitte \`a suspendre l'espace $X$ on peut supposer que tous les produits
de classes de degr\'e strictement positif sont nuls. On peut donc
supposer que~: $$ \alpha_{i,d}^2=0 \  .
$$
Soit $\ell$ tel que $0 \leq \ell \leq d$. D\'efinissons des classes
$\alpha_{i,\ell}$, $i \geq 0$, de degr\'e $2^{k_0+i} +
\ell$, dans $\tilde H^* \Omega ^{d-\ell} X $ comme \'etant les
images, d\'esuspendues $d-\ell$ fois, des classes $\alpha_{i,d}$ par l' homomorphisme de coin
it\'er\'e
$$
Q(H^* X) \ra \Sigma^{d-\ell} H^* \Omega^{d-\ell} X \  .
$$
 Elles sont de degr\'e
 $2^{k_0+i} + \ell$  dans
$\tilde H^* \Omega ^{d-\ell} X $. On va d\'emontrer
qu'elles ont les propri\'et\'es suivantes~:

\begin{itemize}
\item $\alpha_{i,\ell}$ est d\'etect\'ee dans la suite
spectrale d'Eilenberg-Moore de\hfill\break
$\Omega (\Omega^{d-(\ell+1)} X)$ par l'image
des classes $\alpha_{i,\ell+1}$
 dans la colonne $-1$,
\item $Sq^{2^{k_0+i}} \alpha_{i,\ell}=\alpha_{i+1,\ell}$,
\item ces classes sont de degr\'e de nilpotence $\ell$, exactement,
\item $\alpha_{i,\ell}^2 = 0 $, pour tout $i$ assez grand.
\end{itemize}
On ne donne pas de borne pour $i$ dans la
derni\`ere condition.

Pour $\ell=0$, les deux
derni\`eres
 conditions, sont contradictoires. Ceci d\'emontrera le
th\'eor\`eme.

Les propri\'et\'es de ces classes vont \^etre
d\'emontr\'ees par r\'ecurrence descendante.
Les classes $\alpha_{i,d}$ introduites plus haut
satisfont aux conditions requises.
Ceci permet de d\'ebuter la r\'ecurrence.

Supposons donc  les propri\'et\'es \'etablies pour les classes
$\alpha_{i,\ell}, \ell \geq 1$, et d\'emontrons les
pour les classes
$\alpha_{i,\ell-1}$.

\begin{lem}
Les classes $\alpha_{i,\ell -1}$, $i \geq 0$ de
$\bar H^* \Omega^{d-\ell +1} X $ de degr\'e
$2^{k_0+i} + \ell -1$ sont telles que~:
\begin{itemize}
\item $\alpha_{i,\ell -1}$ est d\'etect\'e dans
la suite spectrale d'Eilenberg-Moore de
$\Omega (\Omega ^{d-\ell} X)$ par l'image des classes
$\alpha_{i,\ell}$
 dans la colonne $-1$,
\item $Sq^{2^{k_0+i}} \alpha_{i,\ell -1} = \alpha_{i+1,\ell -1}$,
\item ces classes sont de degr\'e de nilpotence $\ell-1$, exactement,
\item $\alpha_{i,\ell -1}^2 = 0$, pour tout $i$  assez grand.
\end{itemize}
\end{lem}

{\bf D\'emonstration~:~}

La premi\`ere affirmation est cons\'equence de la d\'efinition des classes.

La seconde  affirmation, c'est-\`a-dire la description de l'action
des op\'erations de Steenrod, est une cons\'equence directe des
propri\'et\'es de la suite spectrale d'Eilenberg-Moore ([5], [10])
et de la construction des $\alpha_{i,\ell-1}$ \`a partir de
l'homo\-morphisme de coin. Les suspensions it\'er\'ees
$\Sigma^{d-\ell+1}\alpha_{i,\ell-1}$ sont, par construction, image
par l'homomorphisme de coin it\'er\'e ~:
$$
Q(H^* X) \ra \Sigma^{d-(\ell-1)} H^* \Omega^{d-(\ell-1)} X
$$
des classes $\alpha_{i,d}$.  La relation correspondante a lieu
dans la source, pour les classes
$\alpha_{i,d}$.

{\it A priori} les classes   $  \alpha_{i,\ell}$ sont
au moins $\ell$-nilpotentes et d\'eterminent
des classes    $ \bar \alpha_{i,\ell} \in
R_{\ell}(H^* \Omega^{d-\ell} X )$

Par construction les classes
$ \bar \alpha_{i,\ell-1}$ sont
images des classes  $ \bar \alpha_{i,\ell} \in
R_{\ell}(H^* \Omega^{d-\ell} X )$ par le $F$-isomorphisme fort
$R_{\ell}(H^* \Omega^{d-\ell} X ) \cong_F R_{\ell-1}(H^*
 \Omega^{d-\ell+1} X )$.
Les classes   $  \alpha_{i,d}$ \'etant non-nulles
il en est de m\^eme, par r\'ecurrence, des classes
 $  \alpha_{i,\ell-1}$.

Compte tenu de l'action des op\'erations de Steenrod
l'affirmation concernant le degr\'e de nilpotence
de $\alpha_{i,\ell-1}$ en
r\'esulte.

Il reste \`a montrer  que le cup-carr\'e de ces classes est nul si
$i$ est assez grand. La m\'ethode ne donne pas de borne pour $i$.
C'est, comme dans [8], la partie la plus d\'elicate de la
d\'emonstration,
 l'argument se simplifie
quelque peu, dans la mesure o\`u on ne cherche pas
\`a d\'emontrer de r\'esultat
sur des complexes
finis. Il sera expos\'e dans la prochaine section,
 avec
la partie finale de la d\'emonstration du th\'eor\`eme.

\section {Fin de la d\'emonstration
de  4.1, d\'emonstration  du th\'eor\`eme 0.1}

Pour achever la d\'emonstration on va introduire une famille de
classes $\omega_{i,\ell-1}$ dans la cohomologie de
$\Omega^{d-\ell+1}X$, puis on
\'etudiera
l'action des op\'erations de Steenrod sur ces classes. Ceci
donnera la nullit\'e du cup-carr\'e. Si $\ell=1$, on aura
la contradiction annonc\'ee. On doit \'etudier tous les cas, car
c'est la nullit\'e du cup-carr\'e $\alpha_{i,\ell}^2$ qui entra\^
ine l'existence de $\omega_{i,\ell-1}$.

Supposons donc le cup-carr\'e $\alpha_{i,\ell}^2$ nul
pour $i \geq i_\ell$. Il faut montrer que $\alpha_{i,\ell-1}^2=0$
l'est aussi, pour $i \geq i_{\ell-1} $, o\`u $i_{\ell-1}$ est assez
grand.

  Comme le cup-carr\'e de $\alpha_{i,\ell}$ est nul la classe
$$\alpha_{i,\ell} \otimes
\alpha_{i,\ell} \in (\bar H^*\Omega^{d-\ell}X)^{\otimes 2}
$$ d\'etermine, pour tout $i \geq  i_\ell$,
 un \'el\'ement dans le terme~:
$$E_2^{-2,*}\cong {\rm Tor}^{-2,*}_
{H^*\Omega^{d-\ell}X}
(\F,\F) $$ de la suite
spectrale d'Eilenberg-Moore
de $\Omega(\Omega^{d-\ell}X)$ en degr\'e $2^{k_0+i+1}+2\ell$.

 Cet \'el\'ement
est non-nul. En effet, la classe
$\alpha_{i,\ell} \otimes
\alpha_{i,\ell} \in (\bar H^*\Omega^{d-\ell}X)^{\otimes 2}$
est exactement $2\ell$-nilpotente car elle r\'eduit
non-trivialement dans
$${\Sigma^{2\ell} R_{2\ell}({\rm Tor}^{-2,*}_
{H^*\Omega^{d-\ell}X}(\F,\F))\!  \subset \Sigma^{2\ell}
 R_{2\ell}((\bar H^*\Omega^{d-\ell}X)^{\otimes 2})
\cong  \Sigma^{2\ell}R_{\ell}(\bar H^*\Omega^{d-\ell}X)^{\otimes 2}\ .}
$$
La premi\`ere
inclusion est cons\'equence
de 1.8, car les bords,
qui proviennent de $(\bar H^*\Omega^{d-\ell}X)^{\otimes 3},
$
sont au moins $3\ell$-nilpotents, et comme
$\ell \geq 1$, on peut appliquer la proposition, l'isomorphisme
vient de ce que $\bar H^* \Omega^{d-\ell}X \in {\cal N}il_d$.

\begin{lem} Ce cycle d\'etermine une classe
non-nulle, not\'ee $\omega_{i,\ell-1}$,
dans la cohomologie de $\Omega^{d-\ell+1}X$.
\end{lem}

L'\'el\'ement consid\'er\'e de la colonne $-2$
ne peut \^etre source
d'une diff\'erentielle non-nulle car la colonne $0$ est triviale.

L'image de la diff\'erentielle $d_t$
 dans $E_t^{-2,*}$,
 est un sous-module au moins $((t+2)\ell-t+1)$-nilpotent.  Comme $\ell
\geq 1$, on a $ (t+2)\ell-t+1 > 2\ell$.
On a donc une suite exacte~:
$$
\{0\} \ra D \ra {\rm Tor}^{-2,*}_
{H^*\Omega^{d-\ell}X}(\F,\F)  \ra E_\infty^{-2,*} \ra \{0\} \  ,
$$
o\`u $D$ est au moins $(4\ell-1)$-nilpotent, et
les classes $\alpha_{i,\ell} \otimes
\alpha_{i,\ell}$  sont donc non-triviales dans
l'aboutissement de la suite spectrale.

Les classes $$\omega_{i,\ell-1} \in \bar H^* \Omega^{d-\ell+1}X$$ sont
non-nulles, et
 de degr\'e $2^{k_0+i+1} + 2\ell -2$.
Elles appartiennent au terme $F_{-2}$
de la filtration  d'Eilenberg-Moore
de $\tilde H^* \Omega^{d-\ell+1}X$.
Le terme $F_{-1}$ de la filtration
n'est pas  nul dans ce degr\'e. Il y a donc
 une ind\'etermination sur la d\'efinition des classes
$\omega_{i,\ell-1}$,
seule leur image dans $F_{-2}/F_{-1}$ est bien d\'etermin\'ee.
 L'ind\'etermination
appartient
au terme $F_{-1}$ de la filtration en degr\'e $2^{k_0+i+1} + 2\ell -2$, et
est donc d\'etect\'ee dans $E_\infty^{-1,*}$ en degr\'e $2^{k_0+i+1}
+ 2\ell -1$.

On va calculer l'action de certaines op\'erations de
Steenrod sur les classes $\omega_{i,\ell-1}$. Ce calcul
sera fait modulo l'ordre de nilpotence $2\ell-1$. En fait
ces calculs vont avoir lieu
dans le complexe~:
$$
R_{2\ell-1}(H^* \Omega^{d-\ell}X) \ra R_{2\ell-2}
(H^* \Omega^{d-\ell+1}X)
\ra (R_\ell(H^* \Omega^{d-\ell}X)^{\otimes 2}
$$
qui est exact au centre au sens donn\'e en 1.11. On
a, \`a cause de 1.10 et 3.2, un diagramme commutatif~:
$$\diagram{
R_{2\ell-1}F_{-1}(\ell) &\kern-1pt\ra\kern-1pt & E &\kern-1pt \ra\kern-1pt &  L &\cr \vfl \cong_F&&
\vfl \cong_F && \vfl j & \cr R_{2\ell-2}F_{-1}(\ell-1) &\kern-1pt\ra\kern-1pt
&R_{2\ell-2}F_{-2} (\ell-1)&\kern-1pt\ra\kern-1pt & R_{2\ell-2}
(F_{-2}(\ell-1)/F_{-1}(\ell-1)) &\cr }
$$
o\`u on note $F_{-i}(h)$ pour $F_{-i}(H^* \Omega^{d-h}X)$,
et $L \subset (R_\ell(F_{-1}(\ell))^{\otimes^2} $, et $j$ est un monomorphisme.

Consid\'erons l'image  de la classe $\sigma^{-2\ell}(\alpha_{i,\ell} \otimes
\alpha_{i,\ell} )$ dans $L \subset (R_\ell(F_{-1}(\ell))^{\otimes^2}
$.
Son image par la fl\`eche verticale est l'image
de la classe $\sigma^{-2\ell}\omega_{i,\ell-1}$
dans $R_{2\ell-2} (F_{-2}$ $(\ell-1)
/F_{-1}(\ell-1))$.

\begin{lem}
La classe
$\Sq^{2^{k_0+i}}(\omega_{i,\ell-1})$  ne d\'epend
pas du choix du rel\`evement $\omega_{i,\ell-1}$.
\end{lem}

Par application it\'er\'ee de 3.3 on
sait que $R_{2\ell-2}(F_{-1})$
 dans ${\cal U}_1$. Comme il est par d\'efinition
r\'eduit il est  fortement
$F$-isomorphe \`a une
somme directe de copies de
$F(1)$. Il en r\'esulte
aussit\^ot
que l'op\'eration $\Sq^{2^{k_0+i}}$
est nulle
en degr\'e $2^{k_0+i+1}$. Ce qui donne le r\'esultat
annonc\'e.

Dans la proposition suivante, les \'el\'ements associ\'es \`a
$\alpha_{i,\ell}$ et $\alpha_{i,\ell}
\otimes \alpha_{i,\ell}$  dans
la suite spectrale d'Eilenberg-Moore pour
$\Omega(\Omega^{d-\ell}X)$ sont indiqu\'es entre crochets.

 \begin{prop}
Dans le terme $E_2$ de la suite spectrale on a la relation
$$Sq^{2^{k_0+i}}([\alpha_{i,\ell} \otimes
\alpha_{i,\ell}])=[\alpha_{i+1,\ell} \otimes \alpha_{i,\ell} +
\alpha_{i,\ell} \otimes \alpha_{i+1,\ell}] \  .
$$ Cet \'el\'ement persiste \`a l'infini en une
classe non nulle.
\end{prop}

\begin{cor} Dans la cohomologie
de $\Omega^{d-\ell+1}X$ on a la relation~:
$$Sq^{2^{k_0+i}} \omega_{i,\ell-1}= \alpha_{i+1,\ell-1}
\cup \alpha_{i,\ell-1}  \qquad {\rm mod} \quad (F_{-1})_{2\ell-1}
\ .$$
\end{cor}

La premi\`ere partie de la proposition
 r\'esulte de la
description de l'action des op\'erations de Steenrod dans la suite
spectrale ([5], [10], [11]) et de la formule de Cartan. Comme
seules $\Sq^0$ et $\Sq^{2^{k_0+i}  }$ ont une action non-nulle sur
$u^{2^{k_0+i}}$ la classe~:
$$\Sq^{2^{k_0+i}} ((\omega \otimes u^{2^{k_0+i}}) \otimes
(\omega \otimes u^{2^{k_0+i}} ))
$$
est somme de~:
$$
(\omega \otimes u^{2^{k_0+i+1}} )\otimes
(\omega \otimes u^{2^{k_0+i}} )+
(\omega \otimes u^{2^{k_0+i}} )\otimes
(\omega \otimes u^{2^{k_0+i+1}} ) \  ,
$$
et de~:
$$(\Sq^{2^{k_0+i}}(\omega) \otimes u^{2^{k_0+i}} )\otimes
(\omega \otimes u^{2^{k_0+i}} )+ (\omega \otimes u^{2^{k_0+i}}) \otimes
(\Sq^{2^{k_0+i}}(\omega) \otimes u^{2^{k_0+i}} ) \  .$$ Ces derniers
termes sont nuls, en effet la condition $2^{k-1}\geq h$
 impos\'ee plus haut implique que $\Sq^{2^{k_0+i}}(\omega)=0$.

Cette classe ne peut \^etre l'image d'une diff\'erentielle pour les
m\^emes raisons que $\alpha_{i,\ell} \otimes \alpha_{i,\ell}$. La
proposition en r\'esulte.

Pour ce qui est du corollaire
il r\'esulte des propri\'et\'es
de la suite spectrale d'Eilenberg-Moore
par rapport
au cup-produits, du lemme,  et de ce qui a \'et\'e
dit plus haut~: la formule
a lieu modulo  $F_{-1}$, et ne d\'epend
pas du choix du r\'el\`evement.

\begin{cor} On a ~:
$$  \Sq^{2^k} \Sq^{2^k} \omega_{i,\ell-1} = \alpha_{i+1,\ell-1}^2
\  ,
$$
modulo des termes de degr\'e de nilpotence sup\'erieur ou \'egal
\`a $2\ell-1$.
\end{cor}

La d\'emonstration est identique \`a ce qui pr\'ec\`ede.

\begin{lem} Les classes
$\omega_{i,\ell-1}$ sont
de degr\'e
de nilpotence exactement $2\ell-2$. Il en est de m\^eme
des classes
$\alpha_{i+1,\ell-1}^2$, si elles sont non-nulles
pour tout $i$.
\end{lem}

En effet, les classes $\omega_{i,\ell-1}$ sont non-nulles et on calcule
facilement que
$$
\Sq^{2^{k_0+i+1}}(\omega_{i,\ell-1})=\Sq_{2\ell-2}(\omega_{i,\ell-1})=\omega_{i+1,\ell-1} \  .
$$
Le m\^eme argument s'applique \`a l'autre cas, car on a~:
$$
\Sq^{2^{k_0+i+1}}(\alpha_{i,\ell-1}^2)=
 \Sq_{2\ell-2}(\alpha_{i,\ell-1}^2)=
\alpha_{i+1,\ell-1}^2 \  .
$$
On va en d\'eduire, par l'absurde, la nullit\'e des cup-carr\'es, et
conclure. En fait, l'argument sera le m\^eme, simplement on
utilisera en plus, si $\ell=1$, la relation de cup-carr\'e des
alg\`ebres instables qui donnera la contradiction.

On va commencer par \'enoncer une formule concernant
les op\'erations de Steenrod.

\begin{lem}{\rm[8]}\qua Pour tout entier $n$
$$
\Sq^{2^n}\, \Sq^{2^n} \in \bar {\cal A}(n-1)
\Sq^{2^n} \bar {\cal A}(n-1) \  ,
$$
o\`u ${\cal A}(n-1)$ est la sous-alg\`ebre engendr\'ee par
$\Sq^1,\ldots,\Sq^{2^{n-1}}$ et $\bar{\cal A}(n-1)$
 est l'id\'eal
des \'el\'ements de degr\'e strictement positif.

\end{lem}

La d\'emonstration est bas\'ee sur les relations d'Adem. On
\'ecrit
la relation d'Adem pour $\Sq^{2^n}\, \Sq^{2^n}$~:
$$
\Sq^{2^n}\, \Sq^{2^n} \, = \, \sum_{t=1}^{t={n-1}} \, \Sq^{2^{n+1}-2^t} \Sq^{2^t}  \  .
$$
 Puis
on montre, par r\'ecurrence sur $h$,
que toute op\'eration $\Sq^{2^n+h}$,
avec $0 < h   < 2^{n}$, est dans
$\bar {\cal A}(n-1)
\Sq^{2^n} \bar {\cal A}(n-1) $.

 En fait on a seulement
besoin d'une relation de la forme~:
$$
\Sq^{2^{k_0+i}}\, \Sq^{2^{k_0+i}} = \sum_j \, \, a_j b_j
\  ,
$$
avec $2^{k_0+i+1} > {\rm deg}(b_j)>2^{k_0+i}$. Ceci r\'esulte du lemme, mais
 peut aussi \^etre montr\'e  en utilisant
l'anti-involution $c$ et la base de Cartan-Serre.

On \'ecrit la
d\'ecomposition de $c(\Sq^{2^n}\, \Sq^{2^n})$ sur la base
des mon\^omes admissibles. On observe qu'un mon\^ome admissible
de degr\'e $2^{n+1}$ commence toujours par une op\'eration
de degr\'e strictement sup\'erieur \`a $2^n$. On
observe aussi, en \'evaluant sur $u^{2^{n+1}}$,
 que $\Sq^{2^{n+1}}$ n'appara\^it pas dans la d\'ecomposition.
Puis on r\'eapplique $c$.
%\vskip 5mm

{\bf Fin de la d\'emonstration~:}
Le module  $R_{2\ell-2}(H^*\Omega^{d-\ell+1} X)$, $\ell > 1$ est
non-trivial (\'eventuellement) seulement
 dans les degr\'es de la forme $2^h$ ou $2^h+2^j$.

Pour ce qui est de la cohomologie de $H^*\Omega^dX$, c'est-\`a-dire
si $\ell=1$,
la m\^eme
observation a lieu, si on se restreint au sous-module
$F_{-2}$. Les relations \'etablies plus haut (corollaires 5.3 et 5.4) ont
lieu, \`a d\'esuspension pr\`es,  dans ce sous-module.
On a pour tout $i > i_\ell$ la relation
$$ \Sq^{2^{k_0+i}} \Sq^{2^{k_0+i}} \sigma^{-2\ell+2} (\overline{ \omega_{i,\ell-1}}) =
\sigma^{-2\ell+2} (\overline { \alpha_{i+1,\ell-1}^2})
\  ,
$$
avec les notations \'evidentes.

D'apr\`es le lemme ci-dessus, on a~:
$$
\Sq^{2^{k_0+i}}\, \Sq^{2^{k_0+i}} = \sum_j \, \, a_j b_j
\  ,
$$
 avec $2^{k_0+i+1} >{\rm deg }(b_j)>2^{k_0+i}$. Pour que les classes
 $$
\sigma^{-2\ell+2} (\overline { b_j(\omega_{i,\ell-1}}))
 $$
 soient
 non-nulles elles doivent \^etre
de degr\'e
 de la forme $2^a$, ou $2^a+2^b$.

 Il faut donc que l'on ait une \'equation
 de la forme  ${\rm deg} (b_j)+2^{k_0+i+1}=2^a$,
 soit de la forme ${\rm deg} (b_j)+2^{k_0+i+1}=2^a+2^b$.
La premi\`ere
 \'equation est impossible car ${\rm deg} (b_j)<2^{k_0+i+1}$.
La seconde implique que ${\rm deg} (b_j)$ est une puissance de $2$, ce
qui est \'egalement impossible car $2^{k_0+i+1}>{\rm deg
}(b_j)>2^{k_0+i}$. En fait on a $\alpha({\rm deg}(
\sigma^{-2\ell+2} (\overline { b_j(\omega_{i,\ell-1})})) \geq 3$.

Les classes $\sigma^{-2\ell+2} (\overline { b_j(\omega_{i,\ell-1})})$
sont donc nulles. Il en r\'esulte que les classes
$$
\Sq^{2^{k_0+i}} \Sq^{2^{k_0+i}} ( \sigma^{-2\ell+2}\overline {\omega_{i,\ell-1}})
\in R_{2\ell-2}(H^*\Omega^{d-\ell+1} X)
$$
sont nulles. Donc que les classes
$$
\Sq^{2^{k_0+i}} \Sq^{2^{k_0+i}} ( {\omega_{i,\ell-1}})
$$
sont de degr\'e de nilpotence strictement sup\'erieur \`a $2\ell-2$.

On en d\'eduit, par 1.4, que la classe $\alpha_{i+1,\ell-1}^2$ est nulle d\`es que
$i$ est assez grand. En effet, on a~:
$$
 \Sq_{2\ell-2}^t(\alpha_{i,\ell-1}^2)=
\alpha_{i+t,\ell-1}^2 \  .
$$
Le terme $F_{-2}$ de la filtration de $H^*\Omega^d X$
quotient\'e par le sous-module des \'el\'ements nilpotents
est nul dans les degr\'es qui ne sont pas de la forme $2^h$
ou $2^h+2^j$. Insistons
sur le fait que ceci n'est
pas vrai pour la cohomologie de $\Omega^d X$, mais
la relation utilis\'ee ci-dessous a lieu dans $F_{-2}$.

On a, dans ce sous-module la relation
$$ \Sq^{2^{k_0+i}} \Sq^{2^{k_0+i}} \omega_{i,0} = \alpha_{i+1,0}^2 = \alpha_{i+2,0} \not = 0
\  .
$$
Cette derni\`ere est, exactement, de degr\'e de nilpotence $0$. On
a donc une contradiction. Ce qui ach\`eve la d\'emonstration du
th\'eor\`eme.

En conclusion, on observera que l'on
a utilis\'e nulle part des propri\'et\'es
des foncteurs $R_s$, $s > 2d-1$.
En fait l'argument doit permettre de d\'emontrer la conjecture 0.1.

%\vskip 1cm

\section {D\'emonstration
de  1.5}

On \'enoncera et d\'emontrera dans cette section la proposition
suivante, mal d\'egag\'ee dans [6]~:

\begin{prop}
Soit $M$ un module instable, soit $h \geq 0$ un entier.
Alors l'ensemble des \'el\'ements $x$ tels que  pour
tout $h \geq k \geq 0$ il existe un entier $k_x$, tel que~:
$$
(\Sq_k)^{k_x}x=0 \  ,
$$
est un sous-module de $M$.
\end{prop}

Cette proposition r\'esulte \'evidemment de~:

\begin{lem}
Soit\hspace{-0.5pt} $M$\hspace{-0.5pt} un module instable, soit $h \geq 0$ un entier.
Alors l'ensemble des \'el\'ements $x$ tels que  pour
tout $0 \leq k \leq h$~:
$$
\Sq_kx=0 \  ,
$$
est un sous-module de $M$.
\end{lem}

Ce lemme est cons\'equence des relations d'Adem.
Notons $M_h$ le sous-espace vectoriel gradu\'e d\'etermin\'e par la
condition du
lemme. Il faut montrer que c'est un module sur l'alg\`ebre de
Steenrod.

Soit donc $x \in M_h$. Supposons avoir d\'emontr\'e que pour
$i \leq t-1$ $\Sq^ix \in M_h$, et
montrons que $\Sq^tx \in M_h$.

On calcule $\Sq^{2t}\Sq_k(x)$ \`a l'aide
des relations d'Adem. On a par d\'efinition~:
$$
\Sq^{2t}\Sq_k(x)=\Sq^{2t}\Sq^{|x|-k}(x) \  .
$$
Si $2t < 2(|x|-k)$
les relations d'Adem donnent~:
$$
\Sq^{2t}\Sq_k(x)=
\sum_{0}^{t} \, \, \varepsilon_i
\Sq^{2t+|x|-k-i}\Sq^{i}(x)\  ,
$$ avec $\varepsilon_i=\pmatrix{|x|-k-i-1 \cr
2t-2i \cr
}$.

Si $2t+|x|-k-i > i+|x|$, soit si $2i < 2t-k$, le terme
$\Sq^{2t+|x|-k-i}\Sq^{i}(x)$
est nul
par instabilit\'e. La somme
ci-dessus se r\'eduit donc \`a~:
$$
\Sq^{2t}\Sq_k(x)=
\sum_{  t-{k \over 2} \leq i \leq t} \, \, \varepsilon_i
\Sq^{2t+|x|-k-i}\Sq^{i}(x)\  .
$$
Le coefficient $\varepsilon_t$
est \'egal \`a  $\pmatrix{|x|-k-t-1 \cr
0 \cr
}$ soit \`a $1$, car  $|x|-k-t-1 \geq 0$
par hypoth\`ese. On a donc~:
$$
\Sq^{2t}\Sq_k(x)= \Sq_k \Sq^t(x)+
\sum_{  t-{k \over 2} \leq i \leq t-1} \, \, \varepsilon_i
\Sq_{k-2(t-i)}\Sq^{i}(x)\  .
$$
On en d\'eduit que
$\Sq_k \Sq^t(x)$ est nul  car
 l'hypoth\`ese de r\'ecurrence implique
que tous les autres termes du membre de droite sont nuls,
 et car celui de gauche l'est par hypoth\`ese.

Si $|x|< t+k $ on peut appliquer
les relations d'Adem \`a
$\Sq_k\Sq^t(x)$ directement. On obtient~:
$$
\Sq_{k}\Sq^t(x)=
\sum_{0}^{{t+|x|-k \over 2}} \, \, \epsilon_i
\Sq^{2t+|x|-k-i}\Sq^{i}(x)\  ,
$$ avec $\epsilon_i=\pmatrix{t-i-1 \cr
t+|x|-k-2i \cr
}$.
Soit
$$
\Sq_{k}\Sq^t(x)=
\sum_{t-{k \over 2}\geq i \geq {t+|x|-k \over 2}} \, \, \epsilon_i
\Sq_{k-2(t-i)}\Sq^{i}(x) \  .
$$
Or $i \leq {t+|x|-k \over 2} <  {t}$, et donc $k-2(t-i)<k$,
on peut donc utiliser
l'hypoth\`ese de r\'ecurrence.

\Addresses\recd

\end{document}